%% file: main.tex
\documentclass[12pt]{amsart}
\usepackage{amssymb}
\usepackage{amsmath}
\usepackage{setspace}
\usepackage{mathtools}
\usepackage{verbatim}
\usepackage{amsthm}
\usepackage{tikz-cd}
\usepackage{hyperref}
\usepackage[french, english]{babel}

\newtheorem{theorem}{Theorem}[section]
\newtheorem{lemma}[theorem]{Lemma}
\newtheorem{corollary}[theorem]{Corollary}

\newtheorem{remark}[theorem]{Remark}

\newtheorem{proposition}[theorem]{Proposition}
\newtheorem{defn}[theorem]{Definition}

\newtheorem{example}[theorem]{Example}
\newtheorem{claim}[theorem]{Claim}

\newtheorem{question}[theorem]{Question}
\newcommand{\bb}[1]{\mathbb{#1}}
\newcommand{\sing}[1]{\operatorname{Sing}  #1}
\newcommand{\supp}[1]{\operatorname{Supp} \, #1}
\newcommand{\exc}[1]{\operatorname{Exc} \, #1}

\newcommand{\Null}[1]{\operatorname{Null} \, #1}
\newcommand{\bs}[1]{\operatorname{Bs} \, #1}
\newcommand{\loc}[1]{\operatorname{loc} \, #1}

\newcommand{\cal}[1]{\mathcal{#1}}
\newcommand{\ff}[1]{\mathfrak{#1}}

\newcommand{\Spec}[1]{\text{Spec}\, #1}

\newcommand{\Pic}[1]{\text{Pic}\, #1}

\makeatletter
\newenvironment{abstracts}{%
  \ifx\maketitle\relax
    \ClassWarning{\@classname}{Abstract should precede
      \protect\maketitle\space in AMS document classes; reported}%
  \fi
  \global\setbox\abstractbox=\vtop \bgroup
    \normalfont\Small
    \list{}{\labelwidth\z@
      \leftmargin3pc \rightmargin\leftmargin
      \listparindent\normalparindent \itemindent\z@
      \parsep\z@ \@plus\p@
      
      \itemsep\medskipamount
    }%
}{%
  \endlist\egroup
  \ifx\@setabstract\relax \@setabstracta \fi
}

\newcommand{\abstractin}[1]{%
  \otherlanguage{#1}%
  \item[\hskip\labelsep\scshape\abstractname.]%
}
\makeatother

\DeclareMathOperator{\Aut}{Aut}
\DeclareMathOperator{\lct}{lct}

\title{On the MMP for rank one foliations on threefolds}
\author{Paolo Cascini and Calum Spicer} 
\subjclass[2010]{14E30, 37F75}

\address{Department of Mathematics, Imperial College London, 180 Queen's Gate, 
London SW7 2AZ, UK} 
\email{p.cascini@imperial.ac.uk}

\address{Department of Mathematics, King's College London, Strand,
London WC2R 2LS, UK}
\email{calum.spicer@kcl.ac.uk}

\begin{document}

\begin{abstracts}

\abstractin{english}
We prove existence of flips for log canonical 
foliated pairs of rank one  on a $\mathbb Q$-factorial projective klt threefold. 
This, in particular, provides a proof of the existence of a minimal model for a 
rank one foliation on a threefold for a wider range of singularities, after McQuillan.
\end{abstracts}

\maketitle
\tableofcontents

\input section1.tex

\input section2.tex

\input section3.tex

\input section4.tex

\input section5.tex

\input section6.tex

\input section7.tex

\bibliography{math.bib}
\bibliographystyle{alpha}
\end{document}

%% file: section1.tex
\section{Introduction}

As in the classical Minimal Model Program, it is expected that every foliation on a complex projective manifold $X$ is either uniruled or 
it admits a minimal model, i.e. a birational contraction $X\dashrightarrow X'$ such that the canonical divisor of the induced foliation $\cal F'$ on $X'$ is nef. For rank one foliations on a complex surface, this is known to
be true thanks to the work of Brunella, McQuillan and Mendes (e.g. see \cite{brunella00,mcq08,mendes00}). 
For foliations of rank two on a threefold, the program was carried out in \cite{spicer20,CS18,SS19}.

In \cite{mcq04}, McQuillan proved  the
existence of minimal models
for foliations by curves. More specifically, he showed that if 
$X$ is a projective variety with quotient singularities
and  $\cal F$ is a rank one foliation on $X$ with log canonical singularities, then $\cal F$ admits a minimal model.

\medskip

The goal of this paper is to explore this result in the case of rank one foliations
on threefolds.  In particular, we are interested in proving a generalisation
of McQuillan's theorem, 
and understanding the 
relationship between the birational geometry of foliations and 
classical birational geometry.

In a forthcoming paper \cite{CS24} we show some applications of our results, such as the base point free theorem, the study of foliations with trivial canonical class, and further developing the relationship between the birational geometry of foliations and classical birational geometry.

\subsection{Statement of main results}

Our first main result is to show that flips exist for log canonical foliated pairs of rank one on a $\mathbb Q$-factorial threefold with klt singularities:

\begin{theorem}[= Theorem \ref{thm_flips_exist2}]
Let $X$ be a $\mathbb Q$-factorial klt projective threefold and let $(\cal F,\Delta)$ be a rank one foliated pair on $X$ with  log canonical  singularities. 
Let $R$ be a $(K_{\cal F}+\Delta)$-negative extremal ray such that 
 $\loc R $ has  dimension one (cf. Section \ref{s_nakamaye}).

Then the flipping contraction $\phi\colon X\to Z$ associated to $R$ and the $(K_{\cal F}+\Delta)$-flip exist. 
\end{theorem}

The theorem in particular implies that  the foliated MMP can be run with
natural assumptions on the singularities of the underlying variety,
as well as allowing for the presence of a boundary divisor:

\begin{theorem}[= Theorem \ref{t_mmexist}]
Let $X$ be a $\bb Q$-factorial 
projective threefold with klt singularities and let $(\cal F, \Delta)$ be a  log canonical foliated pair of rank one on $X$.
Assume that $K_{\cal F}+\Delta$ is pseudo-effective.

Then $(\cal F, \Delta)$ admits a minimal model.
\end{theorem}

Our ideas and proofs are greatly indebted to McQuillan's 
strategies and insights, however ultimately our approach
to the existence of minimal models of foliations is independent from the proof given in 
\cite{mcq04} and is based on techniques from the existence of minimal models
in the case of varieties.

\medskip

Finally, we prove several results which relate the birational
geometry of foliations to classical birational geometry.
For instance, it is a striking feature of the canonical model
theorem for foliation by curves on surfaces that 
the singularities on the underlying surface of the canonical
model are never worse than log canonical
(see
 \cite[Fact I.2.4 and Theorem III.3.2]{mcq08}). 

We were interested if such a bound could be proven on threefolds
without making recourse to a canonical model theorem for foliations
on threefolds, which to our knowledge is unknown.
In this direction we prove the following:

\begin{theorem}[=Theorem \ref{t_isolated_sing_bound}]
Let $X$ be a normal threefold and let $\cal F$ be a rank one foliation on $X$ with canonical singularities. Let $0\in X$ be an isolated
singularity.

Then $X$ has log canonical singularities. 
\end{theorem}

Simple examples show that this result is close to optimal in the sense
that if $0 \in X$ is not an isolated singularity then there is in general
no such bound on the singularities of $X$ (see Example \ref{ex_optimal_sing_bound}).

\subsection{Sketch of the proof}

We briefly explain our approach to the proof of existence of flips.
Let $X$ be a $\mathbb Q$-factorial projective threefold with klt singularities and let $\cal F$ be a foliation with canonical singularities on $X$. We assume for simplicity that $\Delta=0$. 
Let $f\colon X \rightarrow Z$ be a $K_{\cal F}$-negative flipping contraction which contracts
a single curve $C$.  We first note that $C$ is necessarily $\cal F$-invariant (cf. \S \ref{s_invariant}).

Our basic approach is to reduce the $K_{\cal F}$-flip to a $(K_X+D)$-flip
for some well chosen divisor $D$ on $X$.  If $D$ is an arbitrary divisor
then there is no reason to expect any relation between $\cal F$
and the pair $(X, D)$.  However, if every component of $D$ 
is $\cal F$-invariant then much of the geometry
of $(X, D)$ is controlled by $\cal F$.

In particular, in Section \ref{s_singularities} 
we show that if $(X, D)$ is log canonical and $C$ is a log canonical centre
of $(X, D)$ then $(K_X+D)\cdot C<0$. 
Thus, the challenge in producing the $K_{\cal F}$-flip becomes to produce
a very singular $\cal F$-invariant divisor containing $C$.  This divisor
gives us the flexibility to produce a divisor $D$ with the desired properties.
This is achieved in Section \ref{s_formalflip}.  
The idea is to perform a careful analysis of 
the singularities of the induced foliation  $f_*\cal F$ on $Z$ at $f(C)$. Unfortunately,
as in the classical MMP, the divisor  $K_{f_*\cal F}$ is not
$\bb Q$-Cartier and so working directly with $f_*\cal F$ is very difficult.  Rather,
we demonstrate the existence of an auxiliary divisor $E$ on $Z$, which is a foliated version of a 
 complement in the classical MMP and such that $K_{f_*\cal F}+E$ is $\bb Q$-Cartier
and the pair $(f_*\cal F, E)$ has mild singularities.  An analysis
of the pair is much more feasible and in fact we are able to 
show that $(f_*\cal F, E)$ admits a particularly simple normal form which, roughly,
can be given by a vector field of the form $\sum n_ix_i\frac{\partial}{\partial x_i}$
where the $n_i$ are non-negative integers.  Examining this normal form, we are able
to produce a large number of invariant divisors containing $C$ on $X$.

It is worth spending a moment to compare this with McQuillan's approach to the existence of a  flip.
In dimension three, it is possible to show 
that $C \cap \sing \cal F$ consists of a single point $P$
and that if $\partial$ is a vector field defining $\cal F$ near $P$
then $\partial = -t\frac{\partial}{\partial t}+ax\frac{\partial}{\partial x}+by\frac{\partial}{\partial y}$
where $C = \{x = y = 0\}$ and $a, b$ are positive integers.  From this, it is possible to deduce 
that the normal bundle of $C$ splits as $\cal O_C(-a)\oplus \cal O_C(-b)$.  By an inductive
analysis
of $\cal F$ along $C$, we can lift this splitting of the normal bundle to a splitting
on a formal neighbourhood of $C$ in $X$, i.e., $C$ is a complete intersection of two formal divisors.
With this description of the formal neighbourhood of $C$ in hand, it is easy to construct a surgery, which is similar to a flip,
by an explicit procedure consisting of a single weighted blow up followed by a single weighted blow down.

\subsection{Acknowledgements}
The first author is partially supported by a Simons collaboration grant. 
The second author is partially funded by EPSRC. We would like to thank 
Florin Ambro, Federico Bongiorno, Mengchu Li, Jihao Liu, James M\textsuperscript cKernan and  Michael McQuillan
for many useful discussions. We are grateful to the referee for carefully reading the paper
and for several useful suggestions and corrections.

%% file: section2.tex
\section{Preliminary Results}

\subsection{Notations}
We work over the field of complex numbers $\mathbb C$. 

Given a normal variety $X$, we denote by $\Omega^1_X$ its sheaf of K\"ahler differentials and,
by $T_X:=(\Omega^1_X)^*$ its tangent sheaf. 
For any positive integer $p$, we denote $\Omega_X^{[p]}\coloneqq (\Omega_X^p)^{**}$. 
Let $A$ be a $\bb R$-Weil divisor on $X$ and let $D$ be a prime divisor.
We  denote by $\mu_DA$ the coefficient of $D$ in $A$.
A log pair $(X,\Delta)$ is a pair of a normal variety and a $\mathbb Q$-divisor $\Delta$ such that $K_X+\Delta$ is $\mathbb Q$-Cartier.  
We refer to \cite{KM98} for the classical definitions of singularities (e.g., klt, log
canonical) appearing in the minimal model program, except for the fact that in our
definitions we require the pairs to have effective boundaries. In addition, we say that a log pair $(X,\Delta)$ is {\bf sub log canonical}, or sub lc, if  $a(E,X,\Delta)\ge -1$ for any geometric valuation $E$
over $X$.
A {\bf fibration} $f\colon X\to Y$ is a surjective morphism between normal varieties with connected fibres. 
We refer to \cite[Section 2.6]{CS18} for some of the basic notions, commonly used in the MMP. 

\medskip

A {\bf foliation of rank $r$} on a normal variety $X$ is a rank $r$ coherent subsheaf $\cal F\subset T_X$ such that 
\begin{enumerate}
\item ${\cal F}$ is saturated in $T_X$, and  
\item $\cal F$ is closed under Lie bracket. 
\end{enumerate}
Note that if $r=1$ then (2) is automatically satisfied. By (1), it follows that $T_X/{\mathcal F}$ is torsion free. 
We denote by $\mathcal N^*_{\cal F}:=(T_X/{\mathcal F})^*$ the {\bf conormal sheaf} of $\cal F$. The {\bf normal sheaf} $\cal N_F$ of $\cal F$ is the dual of the conormal sheaf. 
The {\bf canonical divisor} of $\cal F$ is a divisor $K_{\cal F}$ on $X$ such that $\mathcal O_X(-K_{\cal F})\simeq \det T_{\mathcal F}$.  The foliation $\mathcal F$ is said to be {\bf Gorenstein} (resp.  $\mathbb Q$-Gorenstein) if $K_{\mathcal F}$ is a Cartier (resp. $\mathbb Q$-Cartier) divisor. 
More generally, a {\bf rank $r$ foliated  pair} $(\cal F, \Delta)$ is a pair of a foliation $\cal F$ of rank $r$ and a $\bb Q$-divisor $\Delta\ge 0$
such that $K_{\cal F}+\Delta$ is $\bb Q$-Cartier.

Let $X$ be a normal variety and let $\cal F$ be a rank $r$ foliation
on $X$.  We can associate to $\cal F$ a morphism
\[\phi\colon \Omega^{[r]}_X \rightarrow \cal O_X(K_{\cal F})\]
defined by taking the double dual of the  $r$-wedge product of the map $\Omega^{[1]}_X\to \cal F^*$, induced by the inclusion
$\cal F \subset T_X$.  We will call $\phi$ the {\bf Pfaff field} associated to $\cal F$.
Following \cite[Definition 5.4]{Druel19}, we define the twisted Pfaff field
as the induced map
\[\phi'\colon (\Omega^{[r]}_X \otimes \cal O_X(-K_{\cal F}))^{**} \rightarrow \cal O_X\]
and we define the {\bf singular locus} of $\cal F$,   denoted by $\sing \cal F$, to be the cosupport of the image of $\phi'$.
We say that $\cal F$ is {\bf smooth} at a closed point $x\in X$ if $x\notin \sing \cal F$ and we  say that $\cal F$ is a smooth foliation if $\sing \cal F$ is empty.

Let $\sigma\colon Y\dashrightarrow X$ be a dominant map between normal 
varieties and let $\cal F$ be a foliation of rank $r$ on $X$. We denote by $\sigma^{-1}\cal F$ the {\bf induced foliation} on $Y$ (e.g. see \cite[Section 3.2]{Druel19}).  
If $\sigma\colon Y\to X$ is a  morphism then the induced foliation $\sigma^{-1}\cal F$ is called the {\bf pulled back foliation}. If $f\colon X\dashrightarrow X'$ is a birational map, then the induced foliation on $X'$ by $f^{-1}$ is called the {\bf transformed foliation} of $\cal F$ by $f$ and we will denote it by $f_*\cal F$. 
Moreover, if  $q\colon X' \rightarrow X$ is a quasi-\'etale cover
and $\cal F' = q^{-1}\cal F$ then $K_{\cal F'}=q^*K_{\cal F}$ and 
\cite[Proposition 5.13]{Druel19} implies that $\cal F'$ is smooth if and only if
$\cal F$ is.

\subsection{Singularities in the sense of McQuillan}

The definition of foliation singularities used in \cite{mcq04} is slightly different than the notion defined above. 
We recall McQuillan's definition now.

Let $X$ be a normal variety, let $\mathcal F$ be a rank one foliation on $X$ such that 
$K_{\cal F}$ is $\mathbb Q$-Cartier.  Let $x \in X$ be a point and let $U$ be an open neighbourhood of $x$.  
Up to replacing $U$ by a smaller neighbourhood we may find an index one cover $\sigma\colon U' \rightarrow U$ associated to $K_{\mathcal F}$
and such that $\sigma^{-1}\mathcal F$ is generated by a vector field $\partial$.

We say that $\mathcal F$ is {\bf singular in the sense of McQuillan} at $x \in X$ provided
there exists an embedding $U' \rightarrow M$ where $M$ is a smooth variety and a lift $\tilde{\partial}$
of $\partial$ to a vector field on $M$  such that $\tilde{\partial}$ vanishes at $\sigma^{-1}(x)$.
We denote by $\sing^+{\mathcal F}$ the locus of points $x \in X$ where $\mathcal F$ is singular in the sense 
of McQuillan. Note that $\sing^+\cal F$ does not depend on the choice of $U'$ and it is a closed subset of $X$. 

%On a smooth variety it is easy to see that $\sing^+{\mathcal F} = \sing{\mathcal F}$, but in general
%it is unclear if these two notions of singularity agree. 
 We  have the following inclusion of
singular loci:

\begin{lemma}
\label{lem_mcq_druel_sing}
Let $X$ be a normal variety, let $\mathcal F$ be a rank one foliation on $X$ such that 
$K_{\cal F}$ is $\mathbb Q$-Cartier.

Then $\sing \mathcal F \subset \sing^+{\mathcal F}$.
\end{lemma}
\begin{proof}
See \cite[Lemma 4.1]{CS23b}.
\end{proof}

We will show later that the equality holds if $X$ admits klt singularities (cf. Proposition \ref{p_singcomp}).

\subsection{Invariant subvarieties}\label{s_invariant}
Let $X$ be a normal variety,  and let $\partial \in H^0(X,T_X)$ be a vector field. 
We say that an ideal sheaf $J$ of $X$ is {\bf $\partial$-invariant} if $\partial (J)\subset J$. 
Let $S\subset X$ be a subvariety. Then $S$ is said to be {\bf $\partial$-invariant}, or {\bf invariant by} $\partial$ if
the ideal sheaf ${\cal I}_{S}$ of $S$ is $\partial$-invariant.   

Let $\cal F$ be a foliation on $X$. 
 Then  $S$ is said to be  
{\bf $\cal F$-invariant}, or {\bf invariant by }$\cal F$, if, in a neighbourhood $U$
of the generic point of $S$, $T_{\mathcal F}$ is locally free and for 
any section $\partial \in H^0(U,\cal F)$, we have that $S\cap U$ is $\partial$-invariant. 
If  $D\subset X$ is a prime divisor then we define $\epsilon(D) = 1$ if $D$ is not $\cal F$-invariant and $\epsilon(D) = 0$ if it is $\cal F$-invariant. 

\medskip 

We will  need the following version of Riemann-Hurwitz formula for foliations (e.g. see \cite[Lemma 3.4]{Druel19}): 

\begin{proposition}
\label{prop_RH}
Let $\sigma\colon Y \rightarrow X$ be a finite surjective morphism between normal varieties, let $\cal F$ be a foliation 
on $X$ and let $\cal G := \sigma^{-1}\cal F$.

Then we may write
\[
K_{\cal G} = \sigma^*K_{\cal F} + \sum  \epsilon(\sigma(D))(r_D-1)D\]
where the sum runs over all the prime divisors on $Y$ and $r_D$ is the ramification index of $\sigma$ along $D$.
In particular, if every ramified divisor is $\cal G$-invariant then $K_{\cal G} = \sigma^*K_{\cal F}$.
\end{proposition}

\begin{lemma}\label{l_imginv}
Let $X$ be a normal variety and let $\cal F$ be a rank one foliation on $X$.  
Let $p\colon Y \rightarrow X$ be a proper 
morphism
and assume that $K_{\mathcal F}$ is Cartier and  $K_{\cal G} =  p^*K_{\cal F}$
 where $\cal G := p^{-1}\cal F$.

Then the following hold:
\begin{enumerate}
\item If $W \subset Y$ is a $\mathcal G$-invariant subvariety then $p(W)$ is $\mathcal F$-invariant.

\item If $Z \subset X$ is a $\cal F$-invariant subvariety then $p^{-1}(Z)$
is $\cal G$-invariant.
\end{enumerate}
\end{lemma}
\begin{proof}
We may assume that $X$ is affine and that  $T_{\cal F}$ is generated by a vector field $\partial$
which lifts to a vector field $\tilde{\partial}$ on $Y$ which generates $T_{\cal G}$.

We first prove (1). Let $J$ denote the ideal  of $p(W)$ and let $I$ denote the ideal 
sheaf of $W$. In particular, $p_*I$ is the sheaf associated to $J$. 
%By abusing notation slightly, we have that $p_*I = J$.
Let $f \in J$ and notice that $p^*\partial f = \tilde{\partial}(p^*f)$.
Since $W$ is $\cal G$-invariant and $p$ is proper, we have 
\[
p^*\partial f= \tilde{\partial}(p^*f) \in H^0(Y, I)=J. 
\]
Thus, $J$ is $\partial$-invariant and (1) follows. 

We now prove (2). Let $I$ denote the ideal sheaf of $Z$ and let
$f_1, ..., f_k$ be generators of $I$.  Then $p^*f_1,\dots,p^*f_k$ 
are generators of $p^{-1}I\cal O_Y$, the ideal sheaf
of  the scheme-theoretic preimage $p^{-1}(Z)$. 
%Note that the ideal of  the scheme-theoretic preimage $p^{-1}(Z)$ is the radical of $p^{-1}I\cal O_Y$.
Since $\partial(f_i) \in I$ we get that 
\[
\tilde{\partial}(p^*f_i) = p^*\partial f_i \in p^{-1}I\cal O_Y
\]
and so $p^{-1}I\cal O_Y$ is invariant under $\tilde{\partial}$, as required.
\end{proof}

\iffalse
\begin{lemma}
\label{lem_sing_space_invariant}
Let $X$ be a normal variety, let
$\partial$ be a vector field on $X$
and let $W\subset X$ be a $\partial$-invariant subvariety.

Then $\sing W$ is $\partial$-invariant and, in particular, $\sing X$ is $\partial$-invariant.
\end{lemma}
\begin{proof}
This is \cite[Theorem 5]{MR0212027}
\end{proof}
\fi

\iffalse
\begin{proof}
When $X$ is smooth this is \cite[Lemma 5.1]{CouPer06}. 
We can reduce the general case to the smooth case as
follows.  Everything is local on $X$ so we may freely 
find an embedding $\iota\colon X \hookrightarrow \bb C^M$
and a lift of $\partial$ to a vector field $\tilde{\partial}$ on $\bb C^M$.
Observe that $\iota(W)$ is $\tilde{\partial}$-invariant 
and we may apply the smooth case to conclude.
\end{proof}
\fi

\subsection{Foliation singularities}
Let $X$ be a normal variety and let $(\cal F,\Delta)$ be a foliated pair on $X$.

Given a birational morphism $\pi\colon  \widetilde{X} \rightarrow X$,
let $\widetilde{\cal F}$ be the pulled back foliation on $\tilde{X}$ and let $\Delta'$ be the strict transform of $\Delta$ in $\widetilde X$. 
We may write
\[
K_{\widetilde{\cal F}}+\Delta'=
\pi^*(K_{\cal F}+\Delta)+ \sum a(E, \cal F, \Delta)E\]
where the sum runs over all the prime $\pi$-exceptional divisors of $\tilde X$. 

The rational number $a(E,\cal F,\Delta)$ denotes the {\bf discrepancy} of $(\cal F,\Delta)$ with respect to $E$. 
If $\Delta=0$, then we will simply denote $a(E,\mathcal F)=a(E,\mathcal F,0)$.

\begin{defn}\label{d_canonical} Let $X$ be a normal variety and let $(\cal F,\Delta)$ be a foliated pair on $X$. 
We say that  $(\cal F, \Delta)$ is {\bf terminal} (resp. {\bf canonical},  {\bf log terminal}, {\bf log canonical}) if
$a(E, \cal F, \Delta) >0$ (resp. $\geq 0$, $> -\epsilon(E)$, $\geq -\epsilon(E)$),  for any birational morphism  $\pi\colon \tilde X\to X$ and for any $\pi$-exceptional divisor $E$ on  $\tilde X$.

Moreover, we say that the foliated pair $(\cal F,\Delta)$ is {\bf Kawamata log terminal}, or klt,  if
$\lfloor\Delta\rfloor=0$ and 
if $a(E,\cal F,\Delta)> -\epsilon(E)$ for any birational morphism  $\pi\colon \tilde X\to X$ and for any $\pi$-exceptional prime divisor $E$ on  $\tilde X$.

We say that a $\mathbb Q$-Gorenstein foliation $\cal F$ is terminal (resp. canonical, log canonical) if the foliated pair $(\cal F,0)$ is such. 
\end{defn}

Note that these notions are well defined, i.e., $\epsilon(E)$ and $a(E, \cal F, \Delta)$
are independent of $\pi$.
Observe also that in the case where $\cal F = T_X$, no exceptional divisor $E$ over $X$
is invariant, i.e., $\epsilon(E)=1$ for all $E$, and so this definition recovers the usual
definitions of (log) terminal and  (log) canonical.

Let $P \in X$ be a, not necessarily closed, point of $X$.
We  say  that $(\cal F, \Delta)$ is {\bf terminal} (resp. {\bf canonical, log canonical})
{\bf at $P$} if
for any birational morphism $\pi\colon \widetilde X\to X$ and for any  $\pi$-exceptional  divisor $E$ on $\widetilde X$ whose centre in $X$ is the Zariski closure $\overline P$ of $P$,
we have that the discrepancy of $E$ is $>0$ (resp. $\geq 0$, $\geq -\epsilon(E)$).
Sometimes we will phrase this as $P$ is a terminal (resp. canonical, log canonical) point for $(\cal F, \Delta)$.
We say that $\cal F$ is terminal near $P \in X$ if there is a neighborhood
$U$ of $P$ such that $\cal F\vert_U$ is terminal.
We will see (cf. Lemma \ref{l_terminal}) that being terminal at a closed point $P$ is equivalent to $\cal F$ being
smooth at $P$.

Given an irreducible subvariety $W \subset X$, 
we say that $(\cal F,\Delta)$ is {\bf terminal at the generic point of} 
$W$ if $(\cal F,\Delta)$ is terminal at  the generic point $\eta_W$ of $W$.
We say that $(\cal F,\Delta)$ is {\bf terminal at a general point of} $W$ if $(\cal F,\Delta)$
is terminal  at a general closed point of $W$.

\begin{defn}
Given a normal variety $X$ and a foliated pair $(\cal F, \Delta)$ on $X$, we say that a subvariety $W \subset X$ is a {\bf log canonical
centre} or, in short, {\bf lc centre} (resp. {\bf canonical centre}) of $(\cal F, \Delta)$ if $(\cal F,\Delta)$ is log canonical (resp. canonical) at the generic
point of $W$ and there is a birational morphism $\pi\colon Y\to X$ and a prime divisor $E$ on $Y$ of discrepancy $-\epsilon(E)$ (resp. $0$) whose centre in $X$ is $W$.

A subvariety $W$ is called a {\bf non log canonical centre } of $(\cal F,\Delta)$ if there is a birational morphism $\pi\colon Y\to X$ and a prime divisor $E$ on $Y$ of discrepancy $<-\epsilon(E)$ whose centre in $X$ is $W$.
\end{defn}

Note that if $W$ is a canonical centre of $(\cal F, \Delta)$, then $(\cal F, \Delta)$
is not terminal at the generic point of $W$.
We also remark that if $\cal F$ is smooth and $C\subset X$ is an $\cal F$-invariant curve
then $\cal F$ is terminal at a general point of $C$, but is not terminal
at the generic point of $C$.

\medskip

Given a normal variety $X$ and a foliation $\cal F$ of rank one  on $X$, we say that 
$\cal F$ has {\bf dicritical} singularities if there exists a birational morphism $\pi\colon X'\to X$ and a $\pi$-exceptional divisor $E$ which is not $\pi^{-1}\cal F$-invariant.
We say that $\cal F$ is {\bf non-dicritical}, if it is not dicritical. 

\begin{lemma}\label{l_nondicritical}
Let $X$ be a normal variety and let $\cal F$ be a rank one foliation with canonical singularities.

Then $\cal F$ is non-dicritical.
\end{lemma}
\begin{proof}
This is \cite[Corollary III.i.4]{mp13}.
\end{proof}

Note
that if $\cal F$ is a non-dicritical foliation then the notions of log canonical and canonical coincide. 
In this case
we  might still refer to canonical centres as log canonical centres.
We also remark that any $\cal F$-invariant divisor is an lc centre and a canonical centre 
of $(\cal F, \Delta)$.

\medskip

We will make frequent use of the following consequence of the negativity lemma:

\begin{lemma}\label{l_negativity}
Let $\phi\colon X\dashrightarrow X'$ be a birational map between normal varieties and let 
\begin{center}
\begin{tikzcd}
X \arrow[rd,"f" ' ]\arrow[dashrightarrow]{rr}{\phi} & &X' \arrow{dl}{f'}\\
&Y&  \\
\end{tikzcd}
\end{center}
be a commutative diagram, where $Y$ is a normal variety and $f$ and $f'$ are proper birational morphisms. 
Let $(\cal F,\Delta)$ be a  foliated pair on $X$. Let $\cal F'=\phi_*\cal F$ and let $(\cal F',\Delta')$ be a
foliated pair on $X'$ such that $f_*\Delta=f'_*\Delta'$. 
Assume that $-(K_{\cal F}+\Delta)$ is $f$-ample and $K_{\cal F'}+\Delta'$ is $f'$-ample. 

Then, for any valuation $E$ on $X$, we have
\[
a(E,\cal F,\Delta)\le a(E,\cal F',\Delta').
\]
Moreover, the strict inequality holds if $f$ or $f'$ is not an isomorphism above the generic point of the centre of $E$ in $Y$. 
\end{lemma}
\begin{proof}
The proof is the same as \cite[Lemma 3.38]{KM98}.
\end{proof}

The following is essentially \cite[Corollary III.i.5]{mp13}:

\begin{lemma}\label{l_quasi-etale}
Let $X$ be a normal variety and let $\cal F$ be a rank one foliation on $X$.
Let $q\colon \overline X\to X$ be a finite morphism and let $\overline{\cal F}:=q^{-1}\cal F$. Let $\overline Z\subset \overline X$ be a subvariety and  let $Z:=q(\overline Z)$.  Assume that $(\cal F,\Delta)$ is a foliated pair on $X$  and assume that  $\overline \Delta:=q^*(K_{\cal F}+\Delta)-K_{\overline{\cal F}}$ is an effective $\mathbb Q$-divisor.

Then $(\cal F, \Delta)$ is log canonical at the generic point of $Z$ if and only if 
$(\overline{\cal F}, \overline{\Delta})$ is  log canonical at the generic point of $\overline{Z}$.

Moreover, if $q$ is a quasi-\'etale morphism, then $\cal F$ is terminal (resp. canonical) at the generic point of $Z$ if and only if 
$(\overline{\cal F}, \overline{\Delta})$ is  terminal (resp. canonical) at the generic point of $\overline{Z}$.
\end{lemma}
\begin{proof}
We follow the same methods as \cite[Proposition 5.20]{KM98}. 
Let $\overline f\colon \overline Y\to \overline X$ be a proper birational morphism and let $\overline E$ be an $\overline f$-exceptional divisor on $\overline Y$ whose centre in $\overline X$ is $\overline Z$. Then, by \cite[Theorem VI.1.3]{kollar96}, after possibly replacing $\overline{Y}$ by an higher model, we may assume that 
 there exists a commutative diagram 
 \begin{center}
\begin{tikzcd}
\overline{Y} \arrow[r, "p"] \arrow[d, "\overline f"]
& Y \arrow[d, "f"] \\
\overline X \arrow[r, "q"]
& X
\end{tikzcd}
\end{center}
where $ f$ is birational and $p$ is finite. In particular, if $E=p(\overline{E})$ then $E$ is $f$-exceptional and $Z$ is the centre of $E$ in $X$. 

Assume now that $f\colon Y\to X$ is a proper birational morphism and let $E$ be an $f$-exceptional divisor on $Y$ whose centre in $X$ is $Z$. Let $\overline Y$ be a component of the normalisation of $\overline X\times_X Y$ which maps onto $Y$ and let $\overline {f}\colon \overline Y\to \overline X$  and $p\colon \overline Y\to Y$ be the induced morphisms. Let $\overline E$ be a prime divisor such that $p(\overline E)=E$. 

Lemma \ref{l_imginv} easily implies that $\epsilon(E)=\epsilon(\overline E)$. 
Let $r_E$ be the ramification index of $p$ along $E$. Then, as in the proof of \cite[Proposition 5.20]{KM98}, Proposition \ref{prop_RH} implies that
\[
a(\overline E,\overline{\cal F},\overline {\Delta })=r_Ea(E,\cal F,\Delta)+\epsilon(E)(r_E-1).
\]
It follows easily that 
$a(\overline E,\overline{\cal F},\overline {\Delta }) >  -\epsilon(\overline E)$   
 if and only if $ a(E,\cal F,\Delta)>-\epsilon (E). $ Thus, the first claim follows.

Note that if $q$ is a quasi-\'etale morphism and $\Delta=0$ then $\overline \Delta=0$. 
Lemma  \ref{l_nondicritical} implies that if $\cal F$ (resp. $\overline{\cal F}$) is canonical, then $\epsilon (E)=0$ (resp. $\epsilon(\overline E)=0$). 
Thus, the second claim  follows using the same arguments as above.
 \end{proof}

 Let $f\colon X\to Y$ be a holomorphic morphism between analytic varieties. We say that $f$ is a {\bf submersion} if, for any point $x\in X$, it induces a surjective morphism $df_x:T_xX\to T_{f(x)}Y$. 

\begin{lemma}\label{l_terminal}
Let $X$ be a normal variety and let $\cal F$ be a rank one foliation on $X$ such that 
$K_{\mathcal F}$ is $\mathbb Q$-Cartier. Let $P \in X$ be a closed point.  

Then the following are equivalent:

\begin{enumerate}
\item   $\cal F$ is terminal at $P$. 

\item $P$ is not contained in $\sing^+ \mathcal F$.

\item There is an analytic open neighbourhood $U$ of $P$, 
a quasi-\'etale morphism $q\colon V \rightarrow U$ and a holomorphic submersion $f\colon V \rightarrow B$
such that $q^{-1}\cal F\vert_U$ is induced by $f$.
\end{enumerate}

When $K_{\mathcal F}$ is Cartier these are equivalent to the following:
\begin{enumerate}
\setcounter{enumi}{3}
\item $P$ is not invariant by $\mathcal F$.
\end{enumerate}
\end{lemma}

\begin{proof}
We first observe that all three listed properties are preserved under taking quasi-\'etale covers.
Indeed, terminal singularities are preserved by Lemma \ref{l_quasi-etale}. 
%Being smooth is preserved
%by \cite[Proposition 5.13]{Druel19}.  
Finally, our second and third properties are unchanged by 
a quasi-\'etale cover by definition.

Next, all properties are local about $P$, so we may freely replace $X$ by the index one cover 
associated to $K_{\mathcal F}$ and therefore we may freely assume that $K_{\mathcal F}$ is Cartier.

The equivalence of (2) and (3) is then a consequence of \cite[Lemma I.2.1]{bm16}.

The equivalence of (2) and (4) follows by observing that $P$ is a singular point 
of $\mathcal F$ if and only if $P$ is invariant under $\mathcal F$.

By \cite[Lemma I.1.3]{bm16} if $P$ is invariant, then 
the blow up at $P$ extracts a divisor of discrepancy 
$\leq 0$, in particular $\mathcal F$ is not terminal at $P$. Thus (1) implies (4). 
A direct calculation shows that (3) implies (1). 
\end{proof}

\begin{remark}\label{r_terminal}
Using the same notation as in Lemma \ref{l_terminal}, let $P \in X$ be a point at which $\mathcal F$ is terminal  and let  $C$ be a $\mathcal F$-invariant curve passing through $P$.  Then, for our choice of $q\colon V\to U$, we have that $C':= q^{-1}(C)$ is normal and irreducible and the map $C' \to C$ is ramified over $P$ with ramification index $m$, where $m$ is the Cartier index of $K_{\mathcal F}$.
\end{remark}

Note that the above lemma implies the well known fact
that if $X$ is a surface and $\cal F$ is a terminal rank one foliation
on $X$ then $X$ has at worst quotient singularities.  One can ask more generally if there is a similar way to control 
the singularities of the underlying variety  
in higher dimensions and higher ranks, 
and if such a bound holds if $\cal F$ has only canonical singularities.
For foliations  of co-rank one on a normal threefold, some of these questions were addressed in \cite{CS18}.  
We will approach some cases
of this problem in the rank one case in dimension three (cf. Section \ref{s_singularities}).

We remark that if $\cal F$ is log canonical then there is 
no bound on the singularities of the underlying variety, 
at least from the perspective of Mori theory, as the example in \cite[Example I.2.5]{mcq08} shows.

We also remark that by Lemma \ref{lem_mcq_druel_sing} if $\mathcal F$ is a rank one foliation on a normal variety $X$ such that $\mathcal F$ is terminal at a closed point $P\in X$
then $P \notin \sing \cal F$.

\subsection{Foliations on a surface}
The goal of this section is to present some  results for foliations on a surface which will be used later on. To this end, we employ Mumford's intersection theory for Weil divisors  on a normal projective surface (e.g. see \cite[Example 8.3.11]{fulton84}).

\begin{lemma}
\label{lem_KX+inv_antinef}
Let $X$ be a normal projective surface and let $\cal F$ be a rank one foliation on $X$ such that $K_{\cal F}\equiv 0$ and suppose that $\cal F$ is not algebraically integrable. 
Then
\begin{enumerate}
\item there are only finitely many $\mathcal F$-invariant curves $C_1,\dots,C_k \subset X$; and
\item through a general point of $X$ there exists a curve $M$
not passing through $\sing \cal F$ and such that  $$(K_X+\sum_{i=1}^k C_i)\cdot M \leq 0.$$
\end{enumerate}
\end{lemma}
\begin{proof}
We recall that \cite{jouanolou78} shows that if $X$ is a normal projective surface and $\mathcal F$ is a rank one foliation on $X$ such that $\mathcal F$ is not algebraically integrable, then there are only finitely many $\mathcal F$-invariant curvse on $X$.  This proves item (1). 

We now prove item (2). 
First we show that $\cal F$ has canonical singularities. Suppose not and let $p\colon Y\to X$ be a resolution such that $\cal F_Y:=p^{-1}\cal F$ has canonical singularities, whose existence is guaranteed by Seidenberg's theorem (e.g. see \cite[Theorem 1.1 and pag. 105]{brunella00}).
We have $K_{\cal F_Y}-\sum a(E,\cal F)E\equiv 0$, where the sum runs over all the $p$-exceptional divisors and, by assumption,   there exists a 
$p$-exceptional divisor $E$ such that  $a(E,\cal F)<0$. In particular, $K_{\cal F_Y}$ is not pseudo-effective and by Miyaoka's theorem (e.g. see \cite[Theorem 7.1]{brunella00}), $\cal F_Y$ is  algebraically integrable, and so is $\cal F$, a contradiction. 

Next, observe that we may freely contract $\cal F$-invariant divisors and replace $X$ by a quasi-\'etale cover.
Thus, we are free to assume that $\cal F$ is one of the foliations appearing in
the list \cite[Theorem IV.3.6]{mcq08}.  In particular, $X$ is obtained as an equivariant compactification of a commutative algebraic group of dimension two and $\cal F$ is induced by a codimension one Lie subalgebra.  
We now check each individual case:

\begin{enumerate}
\item $X$ is an abelian surface and $\cal F$ is a linear foliation.  
In particular, if $\cal F$ is not algebraically integrable, there are no $\cal F$-invariant curves on $X$ and $K_X\sim 0$.

\item $X$ is a $\bb P^1$-bundle over an elliptic curve, with projection $p\colon X \rightarrow S$.
In this case, the $\cal F$-invariant curves are either a single section or two disjoint sections. Thus, it is enough to choose $M$ as a general fibre of $p$. 

\item $X$ is a $\bb P^1$-bundle over $\bb P^1$, with projection  $p\colon X\to \bb P^1$.  In this case, the 
$\cal F$-invariant curves are two vertical fibres and either a single or two disjoint sections. Again, we can choose $M$ as a general fibre
of $p$. \qedhere
\end{enumerate}
\end{proof}

\begin{lemma}
\label{lem_KX+inv_antinef2}
Let $X$ be a normal projective surface and let $\cal F$ be a rank one foliation on $X$
which is algebraically integrable.
Let $\Delta, \Theta \geq 0$ be  $\mathbb Q$-divisors on $X$ such that 
\begin{enumerate}
\item $\mu_C \Theta \leq \mu_C\Delta$ for any curve $C$ which is not $\cal F$-invariant, and
\item $(X, \Theta)$ is log canonical.
\end{enumerate}

Then 
$X$ is covered by $\cal F$-invariant curves $M$ such that 
\[
(K_X+\Theta)\cdot M \leq (K_{\cal F}+\Delta)\cdot M.
\]
\end{lemma}

\begin{proof}
%We distinguish two cases. We first assume that $\cal F$ is non-dicritical. Then $\cal F$ is induced by a fibration
%$p\colon X \rightarrow B$. Let   $M$ be a general fibre of $p$. Then 
% $K_X\cdot M = K_{\cal F}\cdot M$ and (1) implies that $\Theta\cdot M \leq \Delta\cdot M$. Thus, our claim follows.
%
%Let us assume now that $\cal F$ is dicritical. 
We may assume, without loss of generality, that the coefficients of $\Delta$ are at most one.  Let $p\colon X' \rightarrow X$ be an F-dlt modification
of $(\cal F,\Delta)$ (cf.  \cite[Theorem 1.4]{CS18}). Then we may write $K_{\cal F'}+p_*^{-1}\Delta+E = p^*(K_{\cal F}+\Delta)$
and $K_{X'}+p_*^{-1}\Theta+E' = p^*(K_X+\Theta)$, where $E, E'$ are $p$-exceptional $\mathbb Q$-divisors and the coefficients of $E$ (resp. $E'$) are greater or equal  (resp. less or equal ) to one.
Since $\cal F'$ is algebraically integrable and non-dicritical, it follows that $\cal F'$ is induced by a fibration $\pi\colon X' \rightarrow B$.
Let $F$ be a general fibre of $\pi$ and observe that 
\begin{enumerate}
%\item $F$ is a rational curve,

\item[(i)] $K_{\cal F'}\cdot F = K_{X'}\cdot F$,

\item[(ii)] $p_*^{-1}\Theta \cdot F \leq p_*^{-1}\Delta \cdot F$, and

\item[(iii)] $E-E' \ge 0$.
%\item $(K_{\cal F'}+p_*^{-1}\Delta+ E)\cdot F = 0$.
\end{enumerate}

Thus, if $M=p(F)$ then 
\[
\begin{aligned}
(K_X+\Theta)\cdot M &=(K_{X'}+p_*^{-1}\Theta+E')\cdot F\\
 &\le (K_{\cal F'}+p_*^{-1}\Delta+E)\cdot F = (K_{\cal F}+\Delta)\cdot M
\end{aligned}
\] and the claim follows.
\end{proof}

\subsection{Adjunction}

\begin{proposition}
\label{prop_adjunction}
Let $X$ be a normal variety and $\cal F$ be a rank one $\bb Q$-Gorenstein foliation on $X$.
Let $S \subset X$
be an $\cal F$-invariant subvariety which is not contained in $\sing \cal F$.  
Let $\nu\colon S^\nu \rightarrow S$
be the normalisation.

Then 
\begin{enumerate}
\item \label{adj} there is an induced foliated pair $(\cal G,\Delta)$ of rank one on $S^\nu$ 
such that
$$K_{\cal F}|_{S^\nu} = K_{\cal G}+\Delta;$$

\item \label{sing} if $(\cal G, \Delta)$ is terminal at a closed point $P \in S^\nu$
then $\cal F$ is terminal at $\nu(P)$.
\end{enumerate}

Assume now that $C \subset X$ is a curve whose irreducible components are 
$\mathcal F$-invariant and they are not contained in $\sing \cal F$. 
If $\nu\colon C^\nu \rightarrow C$ is the normalisation
then $K_{\cal F}|_{C^\nu} = K_{C^\nu}+\Delta$, where $\Delta\ge 0$, and 
\begin{enumerate}
  \setcounter{enumi}{2}
\item $\supp\lfloor \Delta \rfloor = \nu^{-1}(\sing \cal F\cap C)$; and

\item if $P \in C$ is a point such that  $\cal F$ is terminal at $\nu(P)$
then $\mu_P\Delta = \frac{r-1}{r}$  where $r$ is the Cartier index of $K_{\cal F}$ at $\nu(P)$.
\end{enumerate}
\end{proposition}
\begin{proof}
(1) and (2) follow from \cite[Proposition-Definition 3.12]{CS23b} and \cite[Remark 3.13]{CS23b}.

Note that, although \cite[Proposition 3.14]{CS23b} is stated only for  codimension one subvarieties,  the same proof work for any $\cal F$-invariant subvariety. Thus, (3) and, by Remark \ref{r_terminal}, (4) hold. 
\end{proof}

We now explain some generalities comparing foliation 
adjunction and classical adjunction on a threefold:

\begin{proposition}
\label{prop_adj_comp}
Let $X$ be a normal threefold and let 
 $\cal F$ be a foliation of rank one on $X$ with canonical singularities. 
 Let $\Gamma\ge 0$ be a $\mathbb Q$-divisor on $X$ with $\mathcal F$-invariant support and let 
 $S\subset X$ be a reduced and irreducible $\cal F$-invariant divisor such that $(X,\Gamma+S)$ is log canonical. Let $\nu\colon S^\nu\to S$ be its normalisation.

We may write
\[
K_{\cal F}|_{S^\nu} = K_{\cal G}+\Delta
\quad\text{and}\quad
(K_X+\Gamma+S)|_{S^\nu} = K_{S^\nu}+\Theta
\]
where $\cal G$ is the induced foliation and
 $\Delta, \Theta \geq 0$ are $\mathbb Q$-divisors on $S^{\nu}$. Let $C\subset S^{\nu}$ be a curve. 
 
 Then the following hold:

\begin{enumerate}

\item if $\nu(C)$ is contained in $\sing \cal F$ 
then $\mu_C\Delta \geq 1$ and, in particular, $\mu_C\Delta \geq \mu_C\Theta$;

\item if $\nu(C)$ is not contained 
in $\sing \cal F$ and $C$ is not $\cal G$-invariant (i.e., $\cal F$ is terminal at the generic 
point of $\nu(C)$), then 
$\mu_C\Delta = \mu_C\Theta = \frac{n-1}{n}$ where $n$ is the Cartier index of $K_{\mathcal F}$
at the generic point of $C$.
\end{enumerate}
\end{proposition}

\begin{proof}
Let $C\subset S^{\nu}$ be a curve which is not $\cal G$-invariant and such that $\nu(C)$ is not contained in $\sing \cal F$. Then $\nu(C)$ is not contained in the support of $\Gamma$. 

We may calculate $\mu_C\Delta$ using 
\cite[Proposition 3.14]{CS23b}, and $\mu_C\Theta$ by using \cite{kollar13}.
The result then follows.
\end{proof}

Note that, in the notations above,   
if $C$ is $\cal G$-invariant then there is in general no 
natural relation between $\mu_C\Delta$ and $\mu_C\Theta$, as shown in the following example:
\begin{example}
Let $T$ be a smooth surface and let $C_0$ be a smooth curve. Let $X=T\times C_0$ and let $\cal F$ be the foliation induced by the fibration $p\colon X\to T$. 
Let $D\subset T$ be a curve with high multiplicity at a point $z\in D$ and let $S=D\times C_0\subset X$. Then $S$ is $\cal F$-invariant and if  $C=\{z_0\}\times C_0$, we have that  $\mu_C\Delta = 0$,
however  $\mu_C\Theta$ is
arbitrarily large.
\end{example}

\subsection{Jordan decomposition of a vector field}

We follow the notation of \cite{martinet81}. Let $X\coloneqq \widehat{\bb C^m}$ be the completion of $\bb C^m$ at the origin $0\in X$ and let $\partial$ be a 
vector field on $X$ which leaves $W\coloneqq \{0\}$ invariant.  Let $\ff m$ be the maximal ideal defining $W$
and note that, by the Leibniz rule, the ideal  $\ff m^n$ is $\partial$-invariant for all positive integer $n$. Thus, we get a linear map
\[
\partial_n\colon \ff m/\ff m^{n+1} \rightarrow \ff m/\ff m^{n+1}.
\]
We may  write $\partial_n = \partial_{S, n} +\partial_{N, n}$ as the Jordan  decomposition of $\partial_n$
into its semi-simple and nilpotent parts.
This decomposition respects the exact sequences
\[0 \rightarrow \ff m^n/\ff m^{n+1} \rightarrow \bb C[[X]]/\ff m^{n+1} \rightarrow \bb C[[X]]/\ff m^n \rightarrow 0\]
for each positive integer $n$
and it yields a decomposition $\partial = \partial_S+\partial_N$.

We summarise briefly some of the key properties of this decomposition:

\begin{enumerate}
\item $[\partial_S, \partial_N] = 0$;

\item we may find coordinates $y_1, ..., y_m$ on $\widehat{\bb C^m}$ and $\lambda_1,\dots,\lambda_m\in \mathbb C$ so that
$\partial_S = \sum_i \lambda_i y_i \partial_{y_i}$; and 

\item if $Z \subset \widehat{\bb C^m}$ 
is $\partial$-invariant then $Z$ is both $\partial_S$ and $\partial_N$-invariant.
\end{enumerate}

We briefly explain (3).  Let $I_Z\subset \mathbb C[[X]]$ be the ideal of $Z$ and let $I_{Z, n}$ denote
its restriction to $\ff m/\ff m^{n+1}$, for each positive integer $n$. Then $I_{Z, n} \subset \ff m/\ff m^{n+1}$
is a $\partial_n$-invariant subspace and, in particular, it is both $\partial_{S,n}$ and $\partial_{N,n}$-invariant. Thus, (3)  follows.

More generally, we can define the Jordan decomposition for any vector field $\partial$ on the completion of a variety $X$ at a point $P\in X$. Indeed,  consider an embedding $\iota\colon Z \hookrightarrow \bb C^m$ and a lift $\widetilde{\partial}$
of $\partial$ to a vector field on $\bb C^m$.  We can define $\widetilde{\partial}_S$ and $\widetilde{\partial}_N$
as above.  Then $\widetilde{\partial}_S$ and $\widetilde{\partial}_N$ leave $Z$ invariant
and, therefore, they restrict to vector fields $\partial_S$ and $\partial_N$ on $Z$.
Thus,  $\partial = \partial_S+\partial_N$ and  this decomposition has all the properties of the Jordan
decomposition, as described above.

\subsection{Characterising log canonical vector fields}
Let $X$ be a normal variety and let $\partial$ be  a vector field which defines a foliated pair $(\cal F,D)$ such that $K_{\cal F}+D$ is Cartier. Then we say that $\partial$ is {\bf terminal} (resp. {\bf canonical, log canonical}) if the foliated pair $(\cal F,D)$ is such.

Let $P \in Z$ be a germ of a normal variety and let $\partial \in H^0(Z, T_Z)$ be a vector
field which leaves $P$ invariant. By Lemma \ref{l_terminal},  $\partial$ is singular at $P$.
Let $V \coloneqq \ff m/\ff m^2$ where $\ff m$ is the maximal ideal at $P$
and observe that $\partial$ induces a linear map $\partial_0\colon V \rightarrow V$.
Let $\cal F$ be the foliation defined by $\partial$ so that $\partial$
is a section of ${\cal F}(-D)$ for some  
divisor $D \geq 0$.  We assume that $D$ is reduced.

We recall the following results:

\begin{proposition}
\label{prop_lc_is_non_nilp}
Set up as above.

Then the vector field $\partial$ is log canonical at $P$ if and only if $\partial_0$
is non-nilpotent.
\end{proposition}
\begin{proof}
This is   \cite[Fact I.ii.4]{mp13}.
\end{proof}

\begin{proposition}\label{p_semisimple}
Set up as above.  Suppose in addition that either $\partial$ is log canonical and not canonical, or $D \neq 0$.  

Then, after possibly rescaling and taking a change of coordinates, we have that $\partial$ is semi-simple and 
its eigenvalues are all non-negative integers. 
\end{proposition}
\begin{proof}
This follows from \cite[Fact III.i.3]{mp13}. 
\end{proof}

We will also need the following:

\begin{lemma}
\label{lem_eigenvalue_computation}
Let $\partial$ be a log canonical vector field defined over a neighbourhood of 
$0 \in C \subset \bb C^3$ where $C$ is a smooth curve which is invariant by $\partial$. Suppose the following: 
\begin{enumerate}
\item there exist $f_1, ..., f_q$ with $\partial f_i = \lambda_if_i$
where $\lambda_i$ is a positive rational number; and

\item $C$ is an irreducible
component of the reduced locus of $\{f_1 = ... = f_q = 0\}$.
\end{enumerate}

Then (up to rescaling) the semi-simple part of $\partial$
has eigenvalues $1, -a, -b$ where $a, b \in \mathbb Q_{>0}$.
\end{lemma}
\begin{proof}
We may freely replace $\partial$ by its semi-simple part,
 and so we may assume 
that $\partial$ is semi-simple. 
In suitable coordinates and after possibly  rescaling by a unit, we may write 
\[\partial = -x_1\frac{\partial}{\partial x_1}+a_2x_2\frac{\partial}{\partial x_2}+
a_3x_3\frac{\partial}{\partial x_3}\]
and $C = \{x_2 = x_3 = 0\}$

Fix $i\in \{1,\dots,q\}$. By (2), it follows that $f_i \in (x_2, x_3)$, and  we may  write 
\[f_i = \sum_{k,l,m\ge 0} a_{klm}^ix_1^kx_2^lx_3^m\]
for some $a_{klm}^i\in \mathbb C$ such that $a^i_{k00} = 0$ for all $k\ge 0$.
We have 
\[
\partial f_i = \sum a^i_{klm}(-k+a_2l+a_3m)x_1^kx_2^lx_3^m.
\]
Thus,  (1) implies that  
\[
\lambda_i = -k+a_2l+a_3m
\]
for all non-negative integers $k,l,m$ 
such that $a^i_{klm}\neq 0$. 

If $a^i_{kl0}$ (resp. $a^i_{k0m}$) is non-zero for some $i, k, l$ (resp. $i,k,m$)
it follows immediately that $a_2$ (resp. $a_3$) is a positive
rational number.

Assume that $a^i_{k0m} = 0$
for all $i,k, m$.  Then it follows that 
\[
\{x_2 = 0\} \subset \{f_1 = ... = f_k = 0\}
\]
contradicting  the fact that $\{x_2 = x_3 = 0\}$ is an irreducible
component of the latter scheme.  A similar contradiction holds
if $a^i_{km0} = 0$ for all $i, k, m$.
\end{proof}

\subsection{Canonical bundle formula}
\label{s_cbf}
We recall some  results on the canonical bundle formula which will be used later (see \cite{ambro04} for more details). 

Let $(X,\Delta)$ be a sub log canonical  pair and let $f\colon X\to Y$ be a fibration. Assume that 
the horizontal part $\Delta^h$ of $\Delta$ is effective and that 
there exists a $\mathbb Q$-Cartier $\mathbb Q$-divisor $D$ on $Y$ such that 
$$K_X+\Delta\sim_{\mathbb Q} f^*D.$$
If $P$ is a prime divisor on $Y$, we denote by $\eta_P$ its generic point and we define the {\bf log canonical threshold } of $f^*P$ with respect to $(X,\Delta)$ to be 
$$\lct(X,\Delta;f^*P)\coloneqq\sup\{t\in \mathbb R\mid (X,\Delta+tf^*P) \text{ is sub log canonical over $\eta_P$}\}.$$
Let $b_P\coloneqq 1-\lct(X,\Delta;f^*P)$. Then we define  the {\bf discriminant} of $f$ with respect to $\Delta$ as
$B_Y\coloneqq\sum_P b_P P$, where the sum runs over all the prime divisors $P$ in $Y$. 
Let $r$ be the smallest positive integer such that there exists a rational function $\phi$ on $X$ satisfying  
\[
K_X+\Delta+\frac 1 r (\phi)=f^*D.
\]
Then there exists a   $\mathbb Q$-divisor $M_Y$ such that $$K_X+\Delta+\frac 1 r(\phi)=f^*(K_Y+B_Y+ M_Y).$$
$M_Y$ is called the {\bf moduli part} of $f$ with respect to $\Delta$ . 

\iffalse
We will need the following results:
\begin{lemma} \label{l_msemiample}
Set up as above. Assume that $X$ is a threefold,  $Y$ is a surface and $M_Y$ is $\mathbb Q$-Cartier. 

Then $M_{Y}$ is semi-ample. 
\end{lemma}
\begin{proof}
By \cite[Theorem 8.1]{PS09},
there exists a birational morphism $\alpha\colon Y'\to Y$ such that if $X'$ is the normalisation of the main component of $X\times_Y Y'$ with induced morphisms $f'\colon X'\to Y'$ and $\beta\colon X'\to X$, and we write $\beta^*(K_X+\Delta)=K_{X'}+\Delta'$, then the moduli part $M_{Y'}$ of $f'$ with respect to $\Delta'$ is semi-ample. 
Since  $M_{Y}\sim_{\mathbb Q}\alpha_*M_{Y'}$, it follows that the stable base locus of $M_{Y}$ has codimension at least two. Since $Y$ is a surface and $M_Y$ is $\mathbb Q$-Cartier, Fujita Theorem \cite[Theorem 1.10]{fujita83} implies that $M_Y$ is semi-ample, as claimed. 
\end{proof}
\fi

\begin{lemma}
\label{lem_moduli_isotrivial} Let $(X,\Delta)$ be a two dimensional log canonical pair, let $f\colon X\to Y$ be a fibration onto a curve $Y$ and let $D$ be a $\mathbb Q$-divisor on $Y$ such that $K_X+\Delta\sim_{\mathbb Q}f^*D$. 
Let $y\in Y$ be a closed point and assume that there exists an open neighbourhood $U$  of $y$ such that, if we denote 
\[
X_U\coloneqq f^{-1}(U)\qquad \text{and} \qquad X_u\coloneqq f^{-1}(u)\quad \text{for }u\in U
\] 
then $(X_U,\Delta|_{X_U})$ is log smooth and there exists an isomorphism 
\[
\phi_u\colon X_u \to X_y\qquad\text{such that} \qquad \phi^*_u(\Delta|_{X_y})=\Delta|_{X_u} \qquad\text{for all }u\in U.
\]

Then the moduli part of $f$ with respect to $\Delta$ is trivial, i.e.  $M_Y\sim_{\mathbb Q}0$.
\end{lemma}
\begin{proof}
By \cite[Proposition 8.4.9]{Kollar07},
 we may freely 
perform a base change.
Thus, without loss of generality, we may freely
assume that $X \rightarrow Y$ is semi-stable and $\Delta+f^*P$ is a divisor with simple normal crossing for any prime divisor $P$ on $Y$. 

Let $G$ be the support of $\Delta$. By our hypotheses, after possibly replacing $Y$ by a higher cover, we may  find an open subset $V \subset Y$
so that $X = X_0\times V$ and $G = G_0 \times V$, where $X_0$ is a smooth curve and $G_0\subset X_0$ is a finite set.  Since $M_Y$ only depends on the generic fibre we are
therefore free to assume that $X = X_0 \times Y$ and $G = G_0 \times Y$, in which case
the result is immediate.
\end{proof}

\subsection{A recollection on approximation theorems}
\label{s_approx}

We recall some approximation results proven in \cite[Section 4]{CS18}.

We consider the following set up.  Let $\tilde{X} = \Spec{\tilde{A}}$ 
be an affine variety where $\tilde{A}$
is a henselian local ring with maximal ideal $\ff m$ and let
$W \subset \tilde{X}$ be a closed subscheme, defined by an ideal $\tilde{I}\subset \tilde A$.  
Let $\widehat{X} \coloneqq \Spec \widehat{A}$ 
where $\widehat{A}$ is the completion of $\tilde{A}$ along $\tilde{I}$ and let $\widehat{D}$
be a divisor on $\widehat{X}$. Equivalently, $\widehat{D}$ is given by 
a reflexive sheaf $\widehat{M}$ on $\widehat{X}$ and a choice of a section
$\hat{s} \in \widehat{M}$.

The following is a slight generalisation of Artin-Elkik approximation theorem:

\begin{theorem} Set up as above. 
Let $m$ be a positive integer such that $m\widehat{D}$ is Cartier on $\widehat{X} \setminus  W$.

Then, for all positive integer $n$, there exists a divisor $D^n$ on $\tilde{X}$
such that 
\[
D^n = \widehat{D} \mod \tilde{I}^n\qquad\text{and}\qquad
\cal O_{\tilde{X}}(mD^n) \otimes \widehat{A} \cong \cal O_{\tilde{X}}(m\widehat{D}).
\]
\end{theorem}
\begin{proof}See \cite[Corollary 4.5]{CS18}.
\end{proof}

We will use this theorem under the following additional constraints.
Let $X = \Spec A$ be an affine variety and let $P \in X$ be closed point and
suppose $\tilde{A}$ is the henselisation of $A$ at $P$.

\begin{corollary}
Set up as above.  

Then, for all positive integer $n$,
there exists an \'etale neighbourhood  $\sigma \colon U \rightarrow X$ of $P$
and a divisor $D^n_U$ on $U$ such that $\tau^*D^n_U = D^n$ where $\tau\colon \tilde{X} \rightarrow U$
is the induced morphism. 

 In particular, if $\tilde{I} = I \otimes \tilde{A}$
for some $I \subset A$ then $D^n_U = \widehat{D} \mod I^n$.
\end{corollary}

In our applications here we will always take $W = P$ and so the additional hypotheses of the
corollary are always satisfied.

We also recall the following:

\begin{lemma}\label{l_formalklt}
Set up as above.  Suppose in addition that $(\widehat{X}, \widehat{D})$
is klt (resp. (log) terminal, resp. (log) canonical).  

Then for any sufficiently
large positive integer $n$, we have that $(U, D^n_U)$ is klt (resp. (log) terminal, resp. (log) canonical) in a 
neighborhood of $\sigma^{-1}(P)$.
\end{lemma}

\begin{proof}
See \cite[Lemma 4.8]{CS18}.
\end{proof}

\subsection{Resolution of singularities of threefold vector fields}

We recall the following example from \cite{mp13}.

\begin{example} \cite[Example III.iii.3]{mp13} 
\label{ex_unresolvable} Consider the $\bb Z/2\bb Z$-action on $\bb C^3$ given by $(x, y, z) \mapsto (y, x, -z)$.
Let $X$ denote the quotient of $\bb C^3$ by this action.

Consider the vector field on $\mathbb C^3$ given by 
\[\partial_S := x \frac{\partial}{\partial x} - y \frac{\partial}{\partial y}\]
and 
\[\partial_N := a(xy, z)x\frac{\partial}{\partial x} - a(xy, -z)y\frac{\partial}{\partial y} + c(xy, z)\frac{\partial}{\partial z}\]
where $a,c$ are formal functions in two variables such that $c$ is not a unit and it satisfies $c(xy, z) = c(xy, -z)$.
Let  $\partial := \partial_S+\partial_N$.  Note that $\partial \mapsto -\partial$
under the group action. Thus, $\partial$ induces a foliation $\cal F$ on $X$ with an 
isolated canonical singularity and such that $2K_{\cal F}$ is Cartier, but $K_{\cal F}$ is not Cartier.

By  \cite[Possibility III.iii.3.bis]{mp13}, there does not exist a birational morphism $f\colon Y\to X$ such that the induced foliation
$f^{-1}\cal F$ is both Gorenstein and canonical. Moreover, by \cite[III.iii.3.bis]{mp13}, we  also have that the 
curve $\{x= y = 0\}$ is not algebraic, nor analytically convergent.
\end{example}

\begin{defn}\label{d_simple}
Let $X$ be a normal threefold 
 and let $\cal F$ be a rank one foliation on $X$ with canonical singularities. We  say that $\cal F$ admits a {\bf simple singulary} at $P\in X$ if
 either 
\begin{enumerate}
\item ${\cal F}$ is terminal and no component of $\sing  X$ through $P$ is ${\cal F}$-invariant; or

\item $ {X}$ and ${\cal F}$ are formally isomorphic to the variety and the foliation defined 
in Example \ref{ex_unresolvable} at $P$; or

\item $X$ is smooth at $P$.
\end{enumerate}
\end{defn}

\begin{theorem}\label{t_resolution}
Let $X$ be a normal threefold 
 and let $\cal F$ be a rank one foliation on $X$. 

Then there exists a birational morphism (in fact a sequence of weighted blow ups)
$p\colon \widetilde{X} \rightarrow X$ so that $\tilde{\cal F} \coloneqq p^{-1}\cal F$ 
has simple singularities at all points $P \in  \widetilde{X}$.
\end{theorem}
\begin{proof}
%By \cite{bm97}, there exists a resolution of $X$ obtained by only blowing up  foliation invariant centres (see also \cite[Theorem 3.26]{kollar07b}). 
%
Up to replacing $X$ by a resolution of singularities, we may assume that $X$ is smooth.
We may then apply \cite[III.iii.4]{mp13}.
\end{proof}

\begin{lemma}
\label{l_simple_sing_are_quot}
Let $X$ be a normal threefold 
 and let $\cal F$ be a rank one foliation on $X$. 
Suppose that $\cal F$ admits a simple singularity at $P$.  

Then $X$ has cyclic quotient singularities at $P$.
\end{lemma}
\begin{proof}
If $X$ is smooth then there is nothing to show and
if $P \in X$ is as in Example \ref{ex_unresolvable}, then we are done since $X$ is a $\bb Z/2\bb Z$
quotient singularity.

So suppose that $\cal F$ is terminal at $P$. After possibly replacing $X$ by an analytic neighbourhood of $P$, we may assume that there exists a quasi-\'etale cover  $q\colon X' \rightarrow X$ 
 with a holomorphic submersion $f\colon X' \rightarrow S$ as guaranteed by Lemma \ref{l_terminal}.
Assume by contradiction that $X'$ is not smooth.  
Then $q(\sing X') \subset \sing X$ and $q(\sing X')$ is  $\cal F$-invariant, a contradiction.
It follows that  $X'$ is smooth and so $X$ has at worst a cyclic quotient singularity.
\end{proof}
\begin{lemma}

\label{lem_preimagecurve}
Let $G$ a finite group acting on $\mathbb C^3$ without pseudo-reflections,
 let $X \coloneqq \bb C^3/G$ be a quotient singularity and let $q \colon \bb C^3 \rightarrow X$ be the quotient map.
Let $\cal F$ be a   rank one foliation on $X$ and let $C \subset X$ be a smooth $\cal F$-invariant curve.

Then the following hold:
\begin{enumerate}
\item \label{i_case_term} if $\cal F$ is terminal, then $q^{-1}(C)$ is a smooth irreducible curve;

\item \label{i_case_can} if  $q(0)\in X$ is a foliation singularity as in Example \ref{ex_unresolvable}, 
then $q^{-1}(C)$ is a nodal curve and $C$ is a smooth irreducible curve; and 

\item  if the singularity of $\cal F$ at $q(0)$ is simple, then  there is a surface $D\subset X$
 containing $C$ and such that $D$ is klt at $q(0)$ and 
 if $\cal F$ is terminal (resp. canonical) at $q(0)$ then $(D, C)$ is log terminal (resp. log canonical) at $q(0)$. 
\end{enumerate}

\end{lemma}
\begin{proof}
Let $\cal G := q^{-1}\cal F$ and let $C':=q^{-1}(C)$. Then Lemma \ref{l_imginv} implies that  $C'$ is $\cal G$-invariant.  

 If $\cal F$ is terminal then Lemma \ref{l_quasi-etale} implies that $\cal G$ is a terminal foliation on a smooth variety and, by Lemma \ref{l_terminal}, it is 
smooth.  Since $C'$ is a connected leaf of $\cal G$, it  is therefore smooth and irreducible. Thus, (1) follows. 

We now prove (2). Using the same notation as in Example \ref{ex_unresolvable},
we have that $C'$ is necessarily $\partial_S$-invariant.  It is easy to see that the only $\partial_S$-invariant curves passing through $0 \in \bb C^3$ are $\{x = y = 0\}, \{x = z = 0\}$ and $\{y = z = 0\}$.
As in Example \ref{ex_unresolvable},  the curve $\{x= y = 0\}$ is not algebraic, or not analytically convergent.
Thus, $C'$ is either smooth or $C' = \{x = z = 0\} \cup \{y = z = 0\}$ as required. 
Since $C$ is the quotient of $\{x = z = 0\}\cup \{y = z = 0\}$ by the $\bb Z/2\bb Z$-action
we see that $C$ is a smooth irreducible curve. 
Thus, (2) follows. 

Let $D'\subset \mathbb C^3$ be a general surface containing $C'$ and let $D=q(D')$. 
Note that $D'$ is smooth at $0$ and, therefore, $D$ has klt singularities at $q(0)$. 
By \cite[Proposition 5.20]{KM98},  $(D, C)$ is log terminal (resp. log canonical) if and only if $(D', C')$ is log terminal (resp. log canonical).  Thus, (3) follows.
\end{proof}

\subsection{Nakamaye's theorem and the structure of extremal rays}\label{s_nakamaye}

Let $X$ be a normal projective variety and let $M$ be a $\mathbb Q$-Cartier divisor on $X$. 
We define the exceptional locus of $M$ to be 
\[
\Null M\coloneqq \bigcup_{M|_V\text{ is not big }}V
\]
where the union runs over all the subvarieties $V\subset X$ of positive dimension such that $M|_V$ is not big. We denote by $\bb B(M)$ the {\bf stable base locus } of $M$, 
\[
\bb B(M)\coloneqq \bigcap \bs|mM|
\]
where the intersection runs over all the sufficiently divisible positive integers $m$.
Finally, given a ray $R$ in the cone of curves $\overline{\rm NE}(X)$, we define the {\bf locus of $R$} to be the subset
\[
\loc R:=\bigcup_{[C]\in R}C.
\]

 We recall the following result originally due to Nakamaye, in the case of smooth varieties.

\begin{lemma}
\label{lem_nakamaye}
Let $X$ be a normal projective variety. Let $A$ be an ample $\bb Q$-divisor  and let
 $M$ be a big and nef Cartier divisor on  $X$.
 
Then $\Null M = \bb B(M-\epsilon A)$ for any sufficiently small rational number  $\epsilon>0$.
\end{lemma}
\begin{proof}
See \cite[Theorem 1.4]{birkar17}. 
\end{proof}

\begin{proposition}
\label{p_existence_of_alg_contractions}
Let $X$ be a $\mathbb Q$-factorial normal projective variety.  Let
 $M$ be a big and nef Cartier divisor on  $X$. Let $W = \Null M$ and suppose that $M\vert_W \equiv 0$.
 
 Then there exists
a birational contraction $\phi\colon X\to Z$  to an algebraic spaces $Z$, such that $\phi$ contracts  $W$ to a point  and which is an isomorphism outside $W$. 
\end{proposition}
\begin{proof} Let $A$ be an ample divisor. 
Consider the rational map $\phi\colon X \dashrightarrow \bb P^N$ defined by the linear system $|m(M-\epsilon A)|$ where $\epsilon>0$ is a sufficiently small rational number
and $m$ is a sufficiently divisible and large positive integer and note that $\phi$ is birational onto the closure of its image $Y\subset \bb P^N$.
Let $p\colon \overline X\to X$ and $q\colon \overline X\to Y$ be birational morphisms which resolve the indeterminancy locus of $\phi$. 

By Lemma \ref{lem_nakamaye}, it follows that $p(\exc q) = W$. 
We may write 
\[
p^*(m(M-\epsilon A)) = H +F
\]
where $F \geq 0$ is $q$-exceptional
and $H=q^*L$ for some very ample Cartier divisor $L$ on $Y$.
Since $X$ is $\bb Q$-factorial we may choose $G \geq 0$ to be $p$-exceptional so that $-G$ is $p$-ample.
Choose $\delta >0$ sufficiently small so that $A' \coloneqq p^*(m\epsilon A)-\delta G$ is ample.

We therefore have $F+\delta G = p^*(mM) - A'-H$ and  $F+\delta G$ is a $\bb Q$-Cartier $q$-exceptional divisor. 
Since $p(\exc q) = W$, it follows that  $p^*M$ restricted to $\exc q$ is numerically trivial. Thus, if $k$ is a sufficiently divisible positive integer so that   $k(F+\delta G)$ is a Cartier divisor, then
\[
-k(F+\delta G)\vert_{k(F+\delta G)} \equiv k(A'+H)\vert_{k(F+\delta G)}.
\]
Since ampleness of a line bundle on a scheme is equivalent to ampleness of the line bundle
restricted to the reduction and normalisation, and since $A'+H$ restricted to the reduction and normalisation
of each component of $\exc q$ is ample,
we see that $-k(F+\delta G)\vert_{k(F+\delta G)}$ is ample.

We may therefore apply Artin's Theorem \cite[Theorem 6.2]{Artin70} to produce a morphism of algebraic spaces $\overline{X} \rightarrow Z$
which contracts $F+\delta G$ to a point.  By the rigidity lemma this contraction factors through $\overline{X} \rightarrow X$
giving our desired birational  contraction $\phi\colon X \rightarrow Z$.
\end{proof}

%
%
%\begin{proposition}
%\label{prop_loc_R_closed}
%Let $X$ be a normal projective threefold and let  $R$ be an extremal ray of the cone of curves $\overline{\rm NE}(X)$. 
%
%Then $\loc R$ is Zariski closed. 
%\end{proposition}
%\begin{proof}
%Let $H_R$ be the supporting hyperplane to $R$. In particular, $H_R$ is big and nef. By Lemma \ref{lem_nakamaye}, it follows that $\Null H_R$ is a proper Zariski closed subset of $X$. 
%We have $\loc R\subset \Null H_R$. If $\loc R$ is not Zariski closed, then there exists a surface $S\subset \Null H_R$ which contains infinitely many curves in $R$ but such that $S$
% is not contained in $\loc R$. 
%Thus, $H_R\vert_S$ is a nef divisor on $S$ which intersects infinitely many curves trivially and
%therefore it must intersect some moving curve trivially. It follows that $S\subset \loc R$, a contradiction.
%\end{proof}
%

\subsection{Cone theorem}

The cone theorem for rank one foliations
was initially proven in \cite[Corollary IV.2.1]{bm16} when $\cal F$ is Gorenstein
and in \cite{mcq04} when $\cal F$ is $\bb Q$-Gorenstein. A more general version is proven in 
\cite{CS23b}, which we recall here.

\begin{theorem}
\label{t_cone2}
Let $X$ be a normal projective $\bb Q$-factorial variety and let $(\mathcal F, \Delta)$ be a rank one 
foliated pair on $X$. 

Then there are $\mathcal F$-invariant rational curves $C_1,C_2,\dots$ not contained in $\sing^+ {\cal F}$ 
such that 
\[
0<-(K_{\cal F}+\Delta)\cdot C_i\le 2\dim X\]
and
$$\overline{\rm NE}(X)=\overline{\rm NE}(X)_{K_{\cal F}+\Delta\ge 0}+Z_{-\infty}+
\sum_i \mathbb R_+[C_i]$$
where $Z_{-\infty}\subset \overline{\rm NE}(X)$ is a subset contained in the span of the images of
$\overline{\rm NE}(W) \rightarrow \overline{\rm NE}(X)$ where $W \subset X$
are the non-log canonical centres of $(\cal F, \Delta)$.
\end{theorem}
\begin{proof} See \cite[Theorem 4.8]{CS23b}.
\end{proof}

\begin{remark}\label{r_cone}
Set up as in Theorem \ref{t_cone2}. Assume in addition that $(\cal F,\Delta)$ is log canonical and $R$ is a $(K_{\cal F}+\Delta)$-negative extremal ray such that $\dim \loc R=1$. Let  $C$ be a component of $\loc R$. Then \cite[Lemma 4.7]{CS23b} implies that  $C$ is not contained in $\sing^+ {\cal F}$ and, as in the proof of \cite[Theorem 4.8]{CS23b}, we have that $C$ is $\cal F$-invariant. 
\end{remark}

\subsection{A remark on the different notions of singularity}

The following proposition is not needed in this paper, but we believe it is of independent interest as it clarifies the relation between different notions of foliation singularities appearing in the existing literature.

\begin{proposition}\label{p_singcomp}

Let $X$ be a klt variety and let $\mathcal F$ be a rank one foliation on $X$ such that $K_{\mathcal F}$ is $\mathbb Q$-Cartier.

Then $\sing \mathcal F = \sing^+ \mathcal F$.

\end{proposition}

\begin{proof}
By Lemma \ref{lem_mcq_druel_sing} we have the inclusion $\sing \mathcal F \subset \sing^+ \mathcal F$, so suppose for sake of contradiction that
there exists a closed point $x \in \sing^+ \mathcal F \setminus \sing \mathcal F$.   
 We may freely replace $X$ by a neighbourhood of $x \in X$ and we may also
freely replace $X$ be the index one cover associated to $K_{\mathcal F}$.  Thus, we may assume that $\mathcal F$ is defined by a 
vector field $\partial$.
Since $x \not\in \sing \mathcal F$ the morphism $\Omega^{[1]}_X \rightarrow \mathcal O_X$ induced by pairing with $\partial$ is surjective, and so there 
exists a section 
$\omega \in \Omega^{[1]}_X$ such that $\partial(\omega) =1$.
 Let $p\colon X' \rightarrow X$ be a functorial resolution of $X$, cf. \cite[ Theorems 3.35 and 3.45]{kollar07b}. By \cite[Corollary 4.7]{MR2581247} there exists a vector field 
$\partial'$ on $X'$ such that $p_*\partial' = \partial$.  Since $X$ is klt, \cite[Theorem 1.4]{GKKP11} implies that $\omega' \coloneqq p^*\omega$
is a holomorphic 1-form on $X'$.  Note that we still have $\partial'(\omega') = 1$, in particular, $\partial'$ defines 
a smooth foliation $\mathcal F'$ on $X'$.
 
 Since $x \in \sing^+ \mathcal F$, it follows that $x$ is invariant by $\partial$, 
and so $p^{-1}(x)$ is invariant by $\partial'$. 
Perhaps passing to a higher functorial resolution 
we may assume that $p^{-1}(x)$ is a divisor and that there exists an exceptional Cartier divisor $G$ such that $-G$ is $p$-ample. 
 Since $G$ is supported on $p$-exceptional divisors and the $p$-exceptional locus is $\mathcal F'$-invariant we have a 
partial connection $\nabla\colon \mathcal O_{X'}(G) \rightarrow \mathcal O_{X'}(G) \otimes \mathcal O_{X'}(K_{\mathcal F'})$.
Let $E$ be an irreducible component of $p^{-1}(x)$, and let $\mathcal F'_E$ be the restricted foliation.
We may restrict the partial connection $\nabla$ to a partial connection 
\[
\nabla_E \colon \mathcal O_E(G\vert_E) \rightarrow \mathcal O_E(G\vert_E)\otimes \mathcal O_E(K_{\mathcal F'_E}).
\]
Since $\mathcal F'_E$ is smooth, we may apply Bott vanishing to conclude that $G\vert_E^{\dim E} \equiv 0$, cf. \cite[Proposition 5.1]{MR435456}, which contradicts the 
fact that $-G\vert_E$ is ample.
\end{proof}

 In light of this Proposition we ask the following:
 
 \begin{question}
Let $X$ be a normal variety and let $\mathcal F$ be a rank one foliation on $X$ such that $K_{\mathcal F}$ is $\mathbb Q$-Cartier.
Does $\sing \mathcal F = \sing^+ \mathcal F$?
 
 \end{question}

%% file: section3.tex
\section{Facts about terminal singularities}
The following simple observation is a crucial technical ingredient:

\begin{proposition}
\label{prop_admiss_sing_preserved}
Let $X$ be a normal projective variety and let $\cal F$ be a rank one foliation
on $X$ with canonical singularities.
Let $\phi\colon X \dashrightarrow X^+$ be a step of  a 
$K_{\cal F}$-MMP and let  $\cal F^+$ be the induced foliation on $X^+$. 

Then the following hold:
\begin{enumerate}
\item If $X$ admits only quotient singularities, then $X^+$ also admits at worst quotient singularities.
\item If $X$ is a threefold and  $\cal F$ admits simple singularities 
(cf. Definition \ref{d_simple}), then $\cal F^+$ also only admits simple singularities. 

\end{enumerate}
\end{proposition}
\begin{proof}
Let $Z \subset X^+$ be $\phi(\exc \phi)$ if $\phi$ is a 
divisorial contraction and let it be the flipped locus
when $\phi$ is a flip.
In either case by Lemma \ref{l_negativity} if $E$ is a 
divisor centred in a subvariety of $Z$ then $a(E,\cal F^+) >0$.
Thus, $\cal F^+$ is terminal at all points of $Z$, including any generic point of $Z$.

We first prove (1). Assume that $P\in \sing X^+$ is not a quotient singularity.
In particular, $P\in Z$ and $\cal F$ is terminal at $P$. Let $q\colon V\to U$ be a quasi-\'etale morphism over an analytic open neighbourhood 
$U$ of $P$ such that $q^*K_{\cal F^+}$ is Cartier. 
Then $q(\sing V)$ is non-empty.

By Lemma \ref{l_terminal}, after possibly shrinking $U$, we may assume that there exists a submersion  $f\colon V \rightarrow B$ which induces $q^{-1}\cal F^+|_U$
and $\cal F^+$ is not terminal at any generic point of 
$q(\sing V)$. Thus, $q(\sing V)$ is not contained in $Z$. Let $Q\in V$ such that $q(Q)=P$. 
Since $X^+\setminus Z$ has at worst quotient singularities by assumption this implies that $f(Q) \in B$
is a quotient singularity.
Thus, $V$, and hence $U$, has at worst quotient singularities, and (1) follows. 

We now prove (2). Since  $\cal F^+$ is terminal at all points of $Z$, it follows that no components of $Z$ are $\mathcal F$-invariant, so 
if a component $\Sigma$ of $\sing X^+$ is contained in $Z$ then 
$\Sigma$ is not $\cal F^+$-invariant. Thus, (2) follows. 
\end{proof}

\subsection{A version of Reeb stability}
Our goal is to generalise Reeb stability theorem to foliations defined over singular varieties. 

More specifically, let $X$ be a normal variety and let $\cal F$ be a rank one foliation on $X$ which is terminal at all closed points.
Let $C \subset X$ be a compact $\cal F$-invariant curve and let $\Sigma\subset \sing X $ be the locus where $\cal F$ is not Gorenstein. 
By definition of invariance, the set $\{c_1, ..., c_N\}=C\cap \Sigma$ is finite. 
Let $C^\circ = C \setminus \{c_1, ..., c_N\}$ and let $n_k$ be the Cartier index of $K_{\cal F}$ at $c_k$ for each $k=1,\dots,N$. 
We now define the holonomy of $\cal F$ along $C^\circ$. 

Since $C$ is compact, by Lemma \ref{l_terminal}, we may find open sets $U_1,\dots,U_\ell$ in $X$ such that $C$ is contained in the union  $\cup U_i$ and for each $i=1,\dots,\ell$, there exists a finite morphism $q_i\colon V_i\to U_i$ and a fibration $f_i\colon V_i\to T_i$ such that $\cal F_i\coloneqq q_i^{-1}\cal F$ is the foliation induced by $f_i$, $q_i$ is unramified outside $\Sigma$, and if $c_k\in U_i$ for some $k=1,\dots,N$
 then the ramification index of $q_i$ at $c_k$ is $n_k$. In particular, the pre-image of the curve $C$ in $V_i$ is mapped to a point $z_i\in T_i$. 

Pick distinct $i,j$ such that $U_{i,j}:=U_i\cap U_j$ is not empty and it intersects $C$. After possibly shrinking  $U_i$ or $U_j$, we may assume that $U_{i,j}$ does not intersect $\Sigma$. 
Let $V^j_i:=q_i^{-1}(U_{i,j})$ and let $V_{i,j}=V_i^j\times_{U_{i,j}} V_j^i$. Note that the induced morphism $q_{i,j}\colon V_{i,j}\to U_{i,j}$ is unramified and there exists a morphism $f_{i,j}\colon V_{i,j}\to T_{i,j}$ such the pulled back foliation $\cal F_{i,j}$ on $V_{i,j}$ is induced by $f_{i,j}$. Indeed, $f_{i,j}$ is  the Stein factorisation of the morphism $V_{i,j}\to T_i$. 
Let $\sigma_{i,j}\colon T_{i,j}\to T_i$ be the induced morphism. 
Note that the preimage of $C$ in $V_{i,j}$ is mapped to a point $z_{i,j}\in T_{i,j}$ such that $\sigma_{i,j}(z_{i,j})=z_i$.  
After possibly shrinking  $U_i$ and $U_j$, we may assume that $\sigma_{i,j}$ is surjective. 
It follows that $\sigma_{i,j}$ is \'etale. Thus, after replacing $V_i$ by 
$V_{i}\times_{T_{i}}T_{i,j}$
we may assume that $T_i=T_j$.
After repeating this process, finitely many times, we may assume that $T\coloneqq T_i$ and $z\coloneqq z_i\in T$ do not depend on $i=1,\dots,k$. Note that, by the construction above, the choice of the germ $(T,z)$ is uniquely determined by $\cal F$ and $C$. 

 Pick $c\in C^\circ$. Let $\gamma_1,\dots,\gamma_N$ be loops based at $c$  around $c_1,\dots,c_N$ respectively. The {\bf orbifold fundamental group}
 $\pi(C^\circ,c;n_1,\dots,n_N)$
  of $C^\circ$ with weight $n_k$ at $c_k$ is defined as the quotient of $\pi(C^\circ,c)$ by the normal subgroup generated by  $\gamma_1^{n_1},\dots,\gamma_N^{n_N}$.
  We now want to define the {\bf holonomy map}
 $$\rho\colon \pi(C^\circ,c;n_1,\dots,n_N)\to \Aut(T,z),$$
 where $\Aut(T,z)$ denotes the group of biholomorhpic automorphisms on the germ $(T,z)$. 
 Let $\gamma\colon [0,1]\to C^\circ$ be a continuous path which is contained in $U_i$ for some $i=1,\dots,\ell$. Then, since $q_i\colon V_i\to U_i$ is unramified outside $\Sigma$, there exists a lifting $\tilde \gamma\colon [0,1]\to V_i$ of $\gamma$ in $V_i$. Note that $f_i$ maps the image of $\tilde \gamma$ to the point $z\in T$. Proceeding as in the construction of the classic holonomy map (e.g. see \cite{CN85}), we can define a homomorphism
 $$\rho'\colon \pi(C^0,c)\to \Aut (T,z).$$
 Note that if $c_k\in U_i$ for some $i$ and $k$, then the ramification index of $q_i$ at any point in $q_i^{-1}(c_k)$ is equal to $n_k$. Thus, it follows that $\rho'(\gamma_k^{n_k})$ is the identity automorphism of $T$ for any $k=1,\dots,N$ and, in particular, the holonomy map 
  $$\rho\colon \pi(C^\circ,c;n_1,\dots,n_N)\to \Aut(T,z)$$
 is well defined.  

We are now ready to state our singular version of Reeb stability theorem:
\begin{theorem}\label{t_reeb_stability}
Set up  as above. Assume that the image of the holonomy map $\rho$ is finite.

 Then there exists an analytic open subset $W$ of $X$ containing $C$ such that the leaf $C_t$ of $\cal F$ passing through $t\in W$ is a compact analytic subvariety of $W$. 
\end{theorem}

\begin{proof}
The proof of the Theorem is an easy generalisation of the classical Reeb stability theorem  (e.g. see \cite[Theorem IV.5]{CN85}). 
\end{proof}

As a direct application of Reeb stability theorem, we get the following result (see also \cite[II.d.5]{mcq04}):

\begin{proposition}
\label{prop_reeb_stability}
Let $X$ be a normal variety and let $\cal F$ be a rank one foliation on $X$. Let $C \subset X$
be an $\cal F$-invariant curve and suppose that  $\cal F$
is terminal at every closed point $P\in C$.  Suppose moreover that $K_{\cal F}\cdot C<0$.  

Then 
$C$ moves in a family of $\cal F$-invariant curves covering $X$.
\end{proposition}
\begin{proof}
By definition of invariance, $\cal F$ is Gorenstein at the generic point of $C$.
Let $c_1, ..., c_N\in C$ be the non-Gorenstein points of $\cal F$ and let $n_k$
denote the Cartier index of $K_{\cal F}$ at $c_k$, for $k=1,\dots,N$. Let $C^\circ=C\setminus \{c_1,\dots,c_N\}$. 

It follows from foliation adjunction (cf. Proposition \ref{prop_adjunction}),
that $C$ is a rational curve and 
$$K_{\cal F}\cdot C = -2+\sum_{k=1}^N \frac{n_k-1}{n_k}.$$ 
 In particular, since $K_{\cal F}\cdot C<0$
it follows that the orbifold fundamental group $\pi_1(C^\circ,c;n_1,\dots,n_N)$
is finite.  Thus, Theorem \ref{t_reeb_stability} implies the claim. 
\end{proof}

%% file: section4.tex
\section{Subadjunction result in the presence of a foliation}\label{s_singularities}

Given a log pair $(X,S)$, a minimal log canonical centre $W$ of $(X, S)$ and an ample divisor $A$ on $X$, we may write
$(K_X+S+A)|_W = K_W+\Theta$ for an effective divisor $\Theta \geq 0$. We are interested in this situation in the presence of a foliation
which leaves the components of $S$ invariant.  In this case we are able to get some control on $\Theta$ in terms of the singularities
of the foliation.

\subsection{Dlt modification}

Let $X$ be a normal threefold singularity and let $\cal F$ be a rank one foliation 
on $X$ with canonical singularities.
Let $S_1,\dots,S_k$ be prime $\cal F$-invariant divisors.  Our goal here is to control the singularities
of the pair $(X, S\coloneqq \sum a_iS_i)$, where $a_1,\dots,a_k\in (0,1]\cap \mathbb Q$, in terms of the singularities of $\cal F$.
As the following example shows, a canonical foliation singularity
will in general have worse than quotient singularities on the ambient variety  (in contrast to the surface case):

\begin{example}
Let $X = \{xy-zw = 0\} \subset \bb C^4$ 
 and consider the vector field $\partial = x\partial_x-y\partial_y+z\partial_z-w\partial_w$
on $\bb C^4$.  Note that  $X$ is $\partial$-invariant and so $\partial$ induces a rank one foliation
$\cal F$ on $X$.
We claim that  $\cal F$ has canonical singularities.  
Indeed, $\sing \cal F = \{ 0\}$
and if $\ff m$ is the maximal ideal at $0$ then the induced linear map $\ff m/\ff m^2 \rightarrow \ff m/\ff m^2$ is non-nilpotent, 
and Proposition \ref{prop_lc_is_non_nilp} implies that  it is log canonical.
The eigenvalues of $\partial$ are not all positive rational numbers and  \cite[Fact III.i.3]{mp13}  implies that
$\cal F$  has a 
canonical singularity at $(0, 0, 0, 0)$.\end{example}

\begin{lemma}
\label{lem_crep_dlt}
Let $X$ be a normal variety and let $\cal F$ be a rank one foliation on $X$ with canonical singularities.
Let $(X,\Gamma=\sum a_i S_i)$ be a log pair where $S_1,\dots,S_k$ are irreducible $\cal F$-invariant divisors
and $a_1,\dots,a_k\in (0,1]$.

Then there exists a birational morphism  $\mu\colon \overline{X} \rightarrow X$ of 
$(X, \Gamma)$
 such that 
\begin{enumerate}
\item $K_{\overline{\cal F}} = \mu^*K_{\cal F}+F$
where $\overline{\cal F}$ is the foliation induced on $\overline X$ and $F\ge 0$ is a $\mu$-exceptional divisor whose centre in $X$ is contained in the locus where 
$\mathcal F$ is not Gorenstein; and

\item $(\overline X,\overline{\Gamma}+E)$ is dlt and $K_{\overline{X}}+\overline{\Gamma}+E$ 
is nef over $X$, where $E$ is the sum of all the $\mu$-exceptional divisors
and $\overline{\Gamma}$ is the strict transform of $\Gamma$ in $\overline X$.
\end{enumerate}

We  call the morphism $\mu$ a  {\bf dlt modification} of $(X,\Gamma)$ with respect to $\cal F$. 
\end{lemma}

\begin{proof}
Let $U \subset X$ be the Gorenstein locus of $\cal F$, i.e., the open subset of all  points $P \in X$ such that $K_{\cal F}$
is Cartier in a neighbuorhood of $P$ and so $\mathcal F$ is defined by a vector field $\partial$. In particular, $X\setminus U$, being contained in $\sing X$, has codimension at least two. 
Let $p\colon V \to U$ be a functorial resolution, cf. \cite[ Theorems 3.35 and 3.45]{kollar07b}. By \cite[Corollary 4.7]{MR2581247} there exists a lift of $\partial$ to a vector field $\partial'$ on $U$ and therefore 
 we have that  $K_{\cal F_V} = p^*(K_{\cal F}|_U)+F$.
 Since $\cal F$ admits canonical singularities $F = 0$, i.e., $K_{\cal F_V} = p^*(K_{\cal F}|_U)$.

 Let $Y$ be a normal variety which is a partial compactification of $V$ such that there exists a projective morphism $\pi\colon Y\to X$ which extends $p$. Let $\Gamma_Y=\pi^{-1}_*\Gamma$ and let $G$ be the sum of all the $\pi$-exceptional divisors.   Let $Z\to Y$ be a log resolution of $(Y,\Gamma_Y+G)$,  which is an isomorphism over $V$, and let $\rho\colon Z\to X$ be the induced morphism. In particular, if $\mathcal F_Z$ is the  induced foliation on $Z$, $W=\rho^{-1}(U)$ and $q=\rho|_W\colon W\to U$ is the restriction morphism, then $K_{\mathcal F_Z}|_W=q^*(K_{\mathcal F}|_U)$. 
  
We may construct a  morphism $\mu\colon \overline X\to X$ satisfying (2) as the output of an MMP over $X$ starting from $Z$ (e.g. see \cite[Theorem 1.34]{kollar13}).  Let $\overline {\mathcal F}$ be the foliation induced on $\overline X$. It follows that,  if $\overline U=\mu^{-1}U$, then we have that  $K_{\overline {\mathcal F}}|_{\overline U}=r^*(K_{\mathcal F}|_U)$, where $r=\mu|_{\overline U}\colon \overline U\to U$ is the restriction morphism. Thus, since $\cal F$ has canonical singularities, (1) follows. 
\end{proof}

\begin{theorem}
\label{t_isolated_sing_bound}
Let $X$ be a normal threefold and let $\cal F$ be a rank one foliation on $X$ with canonical singularities. Let $0\in X$ be a closed point and 
let $(X, \Gamma)$ be a log pair where $\Gamma$ has $\mathcal F$-invariant support.
%let $(X,\Gamma = \sum a_iS_i)$ be a log pair where $S_1,\dots, S_k$ are irreducible $\cal F$-invariant 
%divisors and $a_1,\dots,a_k\in (0,1]$. 
Suppose that $K_X$ and $\Gamma$ are $\bb Q$-Cartier and that
$(X, \Gamma)$ is log canonical
away from $0$.  
Suppose moreover that $X$ is klt away from $0$.

Then $(X, \Gamma)$ has log canonical singularities. 
\end{theorem}
\begin{proof}
Observe that our hypotheses are preserved by shrinking $X$ and by taking quasi-\'etale covers. Thus,
we may assume without loss of generality that $K_{\cal F}$ is Cartier.

Suppose for the sake of contradiction that $(X, \Gamma)$ has a worse than log canonical singularity at 0.
We may find $0<\lambda <1$, sufficiently close to $1$ so that $(X, \lambda \Gamma)$
is not log canonical at $0$ and is klt away from $0$. Thus, after replacing $\Gamma$ by $\lambda \Gamma$, we may assume that $(X,\Gamma)$ is klt away from $0$. 

Let $\mu\colon \overline{X} \rightarrow X$ be a dlt modification  
of $(X, \Gamma)$ with respect to $\cal F$, whose existence is guaranteed by Lemma \ref{lem_crep_dlt}. Let $\overline {\cal F}\coloneqq\mu^{-1}\cal F$. Then, since $\cal F$ is Gorenstein, 
we have that $K_{\overline{\cal F}}=\mu^*K_{\cal F}$ and since $(X,\Gamma)$ is klt away from $0$, we have that every  $\mu$-exceptional divisor  is centred in $0$.  Let $E=\sum_{i=1}^q E_i$ be the sum of the $\mu$-exceptional divisors and let 
$\overline{\Gamma}$
be the strict transform of $\Gamma$ in $\overline X$.  Lemma \ref{l_nondicritical} implies that $E$ is $\overline{\cal F}$-invariant. 

By classical adjunction and by Proposition \ref{prop_adjunction}, for each $i=1,\dots,q$, we may write
$$(K_{\overline X}+\overline{\Gamma}+E)|_{E_i}=K_{E_i}+\Theta_i\qquad\text{and}\qquad K_{\overline{\cal F}}\vert_{E_i} = K_{\cal G_i}+\Delta_i
$$
for some $\mathbb Q$-divisors $\Delta_i, \Theta_i\ge 0$ on $E_i$ and  where $\cal G_i$ is the induced foliation on $E_i$. 
In particular, $(E_i,\Theta_i)$ is log canonical for all $i=1,\dots,q$.

We first prove the following:

\begin{claim} For any  $i=1,\dots,q$, the surface 
$E_i$ is covered by curves $M$ such that  $(K_{E_i}+\Theta_i)\cdot M \leq 0$.
\end{claim}
\begin{proof}[Proof of the Claim.]
Note that $K_{\cal G_i}+\Delta_i\equiv 0$. 
Suppose first that $\cal G_i$ is not algebraically integrable.
If $\Delta_i\neq 0$, as in  the proof of Lemma \ref{lem_KX+inv_antinef} it follows that $\cal G_i$ is uniruled, a contradiction. Thus, we may assume that  $\Delta_i = 0$, 
and so, by Proposition \ref{prop_adj_comp},
  $\Theta_i$ only consists of $\cal G_i$-invariant components.
Thus, since $(E_i, \Theta_i)$ is log canonical, we have that $\Theta_i \leq \sum C_j$ where the sum runs over all the $\cal G_i$-invariant divisors, and so we may apply Lemma
\ref{lem_KX+inv_antinef} to conclude.

Now suppose that $\cal G_i$ is algebraically integrable. Again, by Proposition \ref{prop_adj_comp}
and since $K_{\cal G_i}+\Delta_i\equiv 0$,
we may apply Lemma \ref{lem_KX+inv_antinef2} 
to conclude. Thus, the claim follows. 
\end{proof}

Let 
$c\colon \overline{X} \dashrightarrow X_{can}$ be the log canonical model
of $(\overline{X}, \overline{\Gamma}+E)$ over $X$, let $\Gamma_{can}\coloneqq c_*\Gamma$ and let $m\colon X_{can} \rightarrow X$
be the induced morphism. 

By (2) of Lemma \ref{lem_crep_dlt}, we have that 
$K_{\overline{X}}+\overline{\Gamma}+E$ is nef over $X$. Thus, 
the inequality of the Claim is in fact an equality and as such, each such curve is contracted by $c$.
This implies that $X_{can} \rightarrow X$ is a small contraction.
In particular, $m^*(K_X+\Gamma) =K_{X_{can}}+\Gamma_{can}$.
Our result follows, since $(X_{can}, \Gamma_{can})$
has log canonical singularities.
\end{proof}

\begin{example}
\label{ex_optimal_sing_bound}
Observe that the assumption that our singularity is isolated in the above theorem is necessary.
Indeed, let $S$ be any normal surface and let $C$ be a smooth curve
and let $\cal F$ be the foliation on $X\coloneqq S \times C$ induced by the projection onto the first coordinate.
It is straightforward to check that $\sing \cal F =\emptyset$ and so $\cal F$ has canonical singularities by \cite[Lemma 5.9]{Druel19},
and moreover, is terminal at all closed points $x \in X$.
\end{example}

We also need  the following:
\begin{proposition}\label{p_boundingsing}
Let $X$ be a normal threefold and let $\cal F$ be a rank one foliation on $X$ with canonical singularities. Let 
$(X,S:=\sum S_i)$ be a log pair where $S_1,\dots, S_k$ are irreducible $\cal F$-invariant divisors and let $C\subset \sing \cal F$ be a curve. 

Then $(X,S)$ is log canonical at the generic point of $C$. 
\end{proposition}
\begin{proof}
The following proof relies on similar, and at the same time easier, ideas as in Theorem \ref{t_isolated_sing_bound}.  Thus, we only sketch its main steps. 

Observe that our hypotheses are preserved by shrinking $X$ and by taking quasi-\'etale covers. Thus,
we may assume without loss of generality that $K_{\cal F}$ is Cartier.

Let $\mu\colon \overline{X} \rightarrow X$ be a dlt modification  
of $(X, S)$ with respect to $\cal F$, whose existence is guaranteed by Lemma \ref{lem_crep_dlt}. Let $\overline {\cal F}\coloneqq\mu^{-1}\cal F$. Then, since $\cal F$ is Gorenstein, 
we have that $K_{\overline{\cal F}}=\mu^*K_{\cal F}$. After possibly replacing $X$ by a neighbourhood of the generic point of $C$, we may assume that  every  $\mu$-exceptional divisor  is centred in $C$.
 Let $E=\sum_{i=1}^q E_i$ be the sum of the $\mu$-exceptional divisors and let 
$\overline{S}$
be the strict transform of $S$ in $\overline X$.  

By classical adjunction and by  Proposition \ref{prop_adjunction}, for each $i=1,\dots,q$, we may write
$$(K_{\overline X}+\overline{S}+E)|_{E_i}=K_{E_i}+\Theta_i\qquad\text{and}\qquad K_{\overline{\cal F}}\vert_{E_i} = K_{\cal G_i}+\Delta_i
$$
for some $\mathbb Q$-divisors $\Delta_i, \Theta_i\ge 0$ on $E_i$ and  where $\cal G_i$ is the induced foliation on $E_i$. 
In particular, $(E_i,\Theta_i)$ is log canonical, for all $i=1,\dots,q$. 

Fix $i=1,\dots, q$ and consider the induced morphism $p\colon E_i\to C$. Let $\Sigma$ be the general fibre of $p$ and let $\Sigma^{\nu}\to \Sigma$ be its normalisation.  
Since $C\subset \sing \cal F$, it follows that a general closed point of $C$ is $\cal F$-invariant. Thus, Lemma \ref{l_imginv} implies that $\Sigma$ is $\overline {\cal F}$-invariant. By classical adjunction and by Proposition \ref{prop_adjunction}, there exist $\mathbb Q$-divisors $\Gamma_i, \Delta'\ge 0$ on $\Sigma^{\nu}$  such that
\[
(K_{E_i}+\Theta_i)|_{\Sigma^\nu} = K_{\Sigma^\nu}+\Gamma_i \qquad \text{and}\qquad  0\equiv K_{\overline{\cal F}}|_{\Sigma^\nu}=K_{\Sigma^\nu}+\Delta'.
\]
By Proposition \ref{prop_adj_comp}, 
it follows that the support of $\Gamma_i$ is contained in the support of 
$\Delta'$ and since $\Delta'$ is integral, whilst $(\Sigma^{\nu},\Gamma_i)$ is 
log canonical, it follows that $\deg(K_{\Sigma^\nu}+\Gamma_i)\le 0$. 
Thus, our results follows as in the proof of Theorem \ref{t_isolated_sing_bound}.
\end{proof}

Note that it is easy to produce examples of a canonical foliation of rank one on a normal variety 
and a collection of invariant divisors $\sum S_i$ so that $(X, \sum S_i)$
has zero-dimensional non-log canonical singularities, as shown in the following example:
\begin{example}
Let $X=\mathbb C^3$, let $\cal F$ be the foliation defined by the  vector field $x\frac{\partial}{\partial x}-y\frac{\partial}{\partial y}$ and let 
$S=\{x = 0\} + \{y = 0\} +\{xy-z^2 =0\}$. Then the support of $S$ is $\cal F$-invariant and  
the origin $0\in X$ is a non-lc centre for $(X,S)$. Note that it is not isolated: the curves $\{x = z = 0\}$ and $\{y = z = 0\}$ are non-lc centres for $(X,S)$ as well.
\end{example}

\begin{remark}
Theorem  \ref{t_isolated_sing_bound}
implies that if $x \in X$ is an isolated $\mathbb Q$-Gorenstein singularity and $\cal F$ is a rank one  foliation with canonical singularities
then $x \in X$ is a log canonical singularity.  It would be interesting to know
if we could improve this bound. E.g. is $x \in X$ log terminal?

Note that if there is a $\bb Q$-Cartier 
$\cal F$-invariant, possibly formal, divisor $S$ passing through $x$ then $(X, tS)$ is log canonical
for $t>0$ sufficiently small
and so $X$ is log terminal.
\end{remark}

\subsection{Subadjuntion}\label{s_subadjunction}

We  work in the following set up. Let $X$ be a 
$\bb Q$-factorial threefold with klt singularities, let $\cal F$ be a rank one foliation on $X$ and let $\Gamma= \sum a_iS_i$ be a $\mathbb Q$-divisor where $S_1,\dots,S_k$ are 
 $\cal F$-invariant prime divisors and $a_1,\dots,a_k\in (0,1)$.  
Let $C \subset X$ be a $\cal F$-invariant projective curve which  is a log canonical centre of $(X, \Gamma)$
and suppose that there are no one-dimensional non-log canonical centres.  
Suppose moreover that $\cal F$ has canonical singularities and that $\cal F$ is terminal 
at a general point of $C$.
Theorem \ref{t_isolated_sing_bound} implies that $(X, \Gamma)$ is log canonical.

By subadjunction for varieties, cf. \cite[Theorem 8.6.1]{Kollar07}, we may write
\[(K_X+\Gamma)\vert_{C^\nu} = K_{C^\nu} +\Theta\]
where $\nu\colon C^\nu \rightarrow C$ is the normalisation and $\Theta\ge 0$ is a $\mathbb Q$-divisor. 

\begin{theorem}
\label{thm_subadj}
Set up as above.  Then
\begin{enumerate}
\item $(C^\nu, \Theta)$ is log canonical;

\item $\lfloor \Theta \rfloor$ is supported on the pre-image of centres of canonical singularities of $\cal F$;

\item if $\cal F$ is terminal at $\nu(Q) \in C$ for some $Q \in C^\nu$ then $\mu_Q\Theta = \frac{n-1}{n}$
where $n$ is the Cartier index of $K_{\cal F}$ at $\nu(Q)$.
\end{enumerate}

In particular, we have 
\[
(K_X+\Gamma)\cdot C\le K_{\cal F}\cdot C.
\]
\end{theorem}

\begin{proof}Let $p\colon \overline{X} \rightarrow X$ be a dlt modification of $(X, \Gamma)$ 
%with respect to $\cal F$ as in Lemma \ref{lem_crep_dlt}
and let $\overline{\Gamma}$ be the strict transform of $\Gamma$ in $\overline X$.
Since $(X,\Gamma)$ is log canonical, we may write 
\[
K_{\overline{X}}+\overline{\Gamma}+E = p^*(K_X+\Gamma)
\]
where $E$ is the sum of all the prime  exceptional divisors of $p$. 
Lemma \ref{l_nondicritical} implies that $E$ is $\overline{\cal F}$-invariant. 
Since $C$ is a log canonical centre of $(X, \Gamma)$, after possibly going to a higher model we may assume that there exists an irreducible component $E_0$ of $E$ 
dominating $C$.  Set $E_1=E-E_0$. 
By  adjunction we may write 
\[
(K_{\overline{X}}+E+\overline{\Gamma})\vert_{E_0} = K_{E_0}+\Theta_0
\]
where $\Theta_0 \geq 0$. 
%Let $\overline{\cal F}$ be the induced foliation on $\overline X$ and let $\cal G$ be the induced foliation  on $E_0$, as in Proposition \ref{prop_adjunction}. 

Let $f\coloneqq p|_{E_0}\colon E_0 \rightarrow C^\nu$ be the restriction morphism. Then $K_{E_0}+\Theta_0$ is $f$-trivial 
and we may write  $K_{E_0}+\Theta_0 = f^*(K_{C^\nu}+M+B)$ where $M\coloneqq M_{C^\nu}$ is the moduli part of $f$  and $B \coloneqq B_{C^\nu} \geq 0$ is the discrepancy part of $f$, as in Section \ref{s_cbf}. 
In particular, we have $\Theta=M+B$. Note that $M$ depends only on $(X, \Gamma)$ in a neighbourhood of the generic point of $C$.  Moreover, for any $P \in C^\nu$, $\mu_PB$ depends only on the germ of $(X, \Gamma)$ at $\nu(P)$.

Since $(E_0,\Theta_0)$ is dlt, it follows that  $(C^\nu, B)$ is log canonical. Moreover, (3) implies (2). Thus, it is enough to prove:

\begin{enumerate}
\item[(a)] \label{m=0}$M = 0$; 
\item[(b)] \label{b_coeff}  for any  closed point $P \in C$ such that  $\cal F$  is terminal at $P$, if $n$ is the Cartier index of $K_{\cal F}$
at $P$, then $\mu_P \Theta = \frac{n-1}{n}$.
\end{enumerate}

We first prove (a). Since $\cal F$ is Gorenstein at the general point $P\in C$ and the support of $\Gamma$ is $\cal F$-invariant,  by Lemma \ref{l_terminal}
there exists an analytic neighbourhood $U$ of $P$ and an isomorphism 
\[
c\colon U \to S \times \mathbb D
\] where $S$ is an analytic surface and  $\mathbb D \subset \mathbb C$ is a disc such that $\cal F|_U$ is induced by the natural submersion $F\colon U\to S$
 and $\Gamma=F^*\Gamma_S$ for some $\mathbb Q$-divisor  $\Gamma_S\ge 0$ on $S$.
Thus, we may assume that $p^{-1}(U)$ is isomorphic to $\overline S\times \mathbb D$ where $\overline S$ is an analytic surface and that $\Gamma+E=\overline F^*D$ for some $\mathbb Q$-divisor $D$ on $\overline S$, where 
$\overline F\colon p^{-1}(U)\to \overline S$ is the natural morphism. 
It follows that for any two general points $P, Q \in C$ we have an isomorphism 
$(f^{-1}(P), \Theta_0\vert_{f^{-1}(P)}) \cong (f^{-1}(Q), \Theta_0\vert_{f^{-1}(Q)})$.   Lemma \ref{lem_moduli_isotrivial} implies that $M = 0$ and (a) follows.

\medskip

We now prove (b). Let $P\in C$ be a closed point such that $\mathcal F$ is terminal at $P$.
By Lemma \ref{l_terminal} there exists  an analytic neighborhood $U$ of $P$ in $X$ and 
a quasi-\'etale cover $q\colon V \rightarrow U$  such that 
  $q^*K_{\cal F}$ is Cartier
  and a holomorphic submersion $F\colon V \rightarrow B$ which induces $\cal F' = q^{-1}\cal F$.
  
    Let $C' = q^{-1}(C)$ and note that $q_C\coloneqq q|_C\colon C' \rightarrow C$
is ramified to order $n$ at $P' \coloneqq q^{-1}(P)$.  We also have that $C'$ 
is $\cal F'$-invariant. 
%As above, since $\cal F'$ is Gorenstein and terminal around $P'$, 
Since $F$ is a submersion, it follows that 
$C'$ is smooth at $P'$.

We may write 
\[
K_V+\Gamma_V=q^*(K_X+\Gamma). 
\]
Note that $C'$ is a log canonical centre for $(V,\Gamma_V)$ and, therefore, by subadjunction  for varieties, we may also write
\[
(K_V+\Gamma_V)|_{C'}=K_{C'}+\Theta',
\]
so that $K_{C'}+\Theta' = q_C^*(K_C+\Theta)$. 

%$K_{F_0}+\Xi_0= g^*(K_{C'}+\Theta')$. 

% we may freely replace $X$ by the index one cover associated to $K_{\mathcal F}$ and so without loss of generality we may assume that $K_{\mathcal F}$ is Cartier at $P$.  
Since $\Gamma_V$ is $\cal F'$-invariant, after replacing $U$ by a smaller analytic neighbourhood of $P'$, we have that the submersion $F$ defines an isomorphism  
\[
c\colon V \to S \times \mathbb D
\] where $S\subset B$ is an analytic open set,  $\mathbb D \subset \mathbb C$ is a disc and $\Gamma_V=F^*\Gamma_S$ for some $\mathbb Q$-divisor  $\Gamma_S\ge 0$ on $S$. It follows that $\mu_{P'}\Theta'=0$ and, therefore, by Riemann-Hurwitz we have that 
$\mu_P\Theta = \frac {n-1}n$, as claimed. 
This concludes the proof of  (b).  Thus, (1), (2) and (3) follow.

\medskip

Our final claim follows immediately from the results above and  Proposition \ref{prop_adjunction}.
\end{proof}

%% file: section5.tex
\section{The formal neighborhood of a flipping curve}\label{s_formalflip}

Let $X$ be a normal threefold, let $\cal F$ be a rank one foliation on $X$
 and let $f\colon X \rightarrow Z$ be a $K_{\cal F}$-flipping contraction.  
 Let $C = {\rm Exc} (f)$.
In the case where $C$ is smooth and irreducible, McQuillan has produced a 
rather complete picture of the structure
of a formal neighborhood of $C$ by examining the formal holonomy around the curve, 
in particular, he shows the existence of a formal $\cal F$-invariant divisor
containing $C$.

In this section we provide a different approach to producing such an invariant divisor.
Our two main ingredients are a foliated analogue of the existence of complements and an analysis
of the structure of log canonical foliation singularities.  

\subsection{Preliminary results}
We begin with the following results:

\begin{lemma}\label{l_generalelement}
Let $X$ be a normal threefold with only quotient singularities and let 
$C \subset X$ be a curve such that $\sing X \cap C = \{x_1, ..., x_k\}$ is finite. 
Let $H$ be an ample divisor and  assume that for each $i=1,\dots,k$, 
there exists  a prime divisor $D_i$ such that $D_i$ is klt at $x_i$, 
$C$ is contained in $D_i$, and the log pair $(D_i,C)$ is log canonical at 
the point $x_i$. 

Then, after possibly replacing $X$ by an analytic neighbourhood of $C$, there exists a divisor $L$ such that for any sufficiently large positive integer $m$ 
the general element $D$ of the linear system 
\[
\{ \Sigma \in |L+mH|\mid C\subset \Sigma\}
\]
is such that $D$ is klt at each point $x_1,\dots,x_k$ and $(D,C)$ is log canonical. 
\end{lemma}

\begin{proof}
After possibly replacing $X$ by an analytic neighbourhood of $C$, for each $i=1,\dots, k$ we may find an effective divisor $M_i$ on $X$, such that
$$\sum_{j\neq i} D_j + M_i \text{ is Cartier at $x_i$} \qquad \text{ and }\qquad  x_j\notin M_i \text{ for any } j\neq i.$$
Let $L:=\sum_{i=1}^k (D_i+M_i)$. Then $L-D_i$ is Cartier at $x_i$ for each $i=1,\dots,k$. Thus, if $m$ is a sufficiently large positive integer, we have that $x_i$ is not contained in the base locus of $|L-D_i+mH|$. In particular, there exists $\Sigma_i\in |L+mH|$ such that $C\subset \Sigma_i$, $\Sigma_i$ has klt singularities at $x_i$ and $(\Sigma_i,C)$ is log canonical at $x_i$. Thus, the general element in the linear system 
\[
\{\Sigma\in |L+mH|\mid C\subset \Sigma\}
\] satisfies the required properties. 
\end{proof}

\begin{lemma}
\label{lem_pic_inj}
Let $X$ be a  normal variety of dimension at least three and with only quotient singularities and let $C$ be a Cartier divisor on $X$.
Let $H$ be an ample divisor on $X$, let $L$ be a  divisor on $X$
and let $D\in \lvert mH+L\rvert$  for a sufficiently large positive integer $m$. Suppose that $C\vert_D \sim_{\mathbb Q} 0$.  

Then $C \sim_{\mathbb Q} 0$.
\end{lemma}

\begin{proof}
After replacing $C$ by a multiple, we may assume that $C|_D\sim 0$ and that there exists 
a compactification $\overline X$ of $X$ which is normal and it admits a Cartier divisor $\overline C$ such that $\overline C|_X=C$. 

Let $\pi \colon Y \rightarrow \overline{X}$ be a finite cover which is unramified along the general point of $D$ and such 
that $L_Y\coloneqq \pi^*L$ is Cartier. 
Let $C_Y \coloneqq \pi^*\overline{C}$ 
and $D_Y \coloneqq \pi^{-1}(D)$.  Notice that $C_Y\vert_{D_Y} \sim 0$. Let $\overline{D_Y}$ be the
closure of $D_Y$ in $Y$. 
It follows that $C_Y|_{\overline {D_Y}}\sim\sum  a_i C_i|_{\overline{D_Y}}$ where $a_i\in \mathbb Z$ and $C_i$ is a divisor contained in $Y\setminus U$, where $U \coloneqq \pi^{-1}(X) \subset Y$. 

By choosing $m \gg 0$ we may assume by Serre duality and Serre vanishing that 
\[
H^1(Y, \cal O_Y(-\overline{H_Y})) = H^2(Y, \cal O_Y(-\overline{H_Y})) =  0.
\] 
By the exact sequence
$$1 \rightarrow \cal O_Y(-\overline{D_Y}) \rightarrow \cal O^*_{Y} \rightarrow \cal O^*_{\overline{D_Y}} \rightarrow 1$$
it follows that 
$\Pic Y \rightarrow \Pic \overline{D_Y}$ is an isomorphism. 
Thus, $C_Y\sim \sum a_i C_i$ and, in particular, $C_Y\vert_U\sim 0$. 

Perhaps passing to the Galois closure of $U \rightarrow X$ with Galois
group $G$ we see that if $s \in H^0(U, C_Y\vert_U)$ is a non-vanishing section then $\prod_{g\in G} g\cdot s$
is a nonvanishing $G$-invariant section of $q C_Y\vert_U$, where $q=\# G$, and so descends
to a nonvanishing section of $q C$. 
Thus, the claim follows. 
\end{proof}

\subsection{Flipping contractions}
Let $X$ be a projective $\mathbb Q$-factorial normal threefold and 
let $\cal F$ be a rank one foliation on $X$.
% with log canonical singularities.  
Let $R$ be a $K_{\cal F}$-negative extremal ray and assume 
that $\loc R$ is a curve $C$. %By Theorem \ref{t_cone2} and Remark \ref{r_cone}, it follows that every component of $C$ is $\cal F$-invariant. 
Let $H_R$ be a supporting hyperplane to $R$ for $\overline{NE}(X)$.

\begin{lemma}
\label{lem_null=loc}Set up as above. 
Let $S$ be a surface.  

Then $H_R^2\cdot S>0$ and, in particular, $\Null H_R = \loc R$.
\end{lemma}
\begin{proof}
Suppose for the sake of contradiction that $H_R^2\cdot S= 0$.

Let $\nu\colon S^\nu \rightarrow S$ be the normalisation of $S$.
Since $H_R$ is big and nef we may write $H_R \sim_{\bb Q} A+B+tS$ where $A$ is ample, $B \geq 0$ and does not contain $S$
in its support and $t>0$. 
It follows
that 
\[
\nu^*H_R\cdot \nu^*S=\frac 1 t \nu^*H_R\cdot \nu^*(H_R-A-B)
<0.
\]

We may also write $H_R \sim_{\bb Q} K_{\cal F}+A'$ where $A'$ is ample. Since $H_R^2\cdot S = 0$
we see that $\nu^*H_R\cdot \nu^* K_{\cal F} = -\nu^*H_R\cdot \nu^*A' <0$.

Suppose first that $S$ is not $\cal F$-invariant. Then 
\cite[Proposition-Definition 3.7]{CS23b} implies that there exists a $\mathbb Q$-divisor $D\ge 0$ on $S^{\nu}$ such that 
$(K_{\cal F}+S)\vert_{S^\nu} \sim_{\mathbb Q} D $.
We have
\[0 \leq \nu^*H_R\cdot \nu^*(K_{\cal F}+S)  = (-\nu^*H_R\cdot \nu^*A') + (\nu^*H_R\cdot \nu^*S) <0, \]
which gives us a contradiction.

Thus,  we may assume that $S$ is $\cal F$-invariant. Let $(\cal G,\Delta)$ be the induced foliated pair on $S^\nu$, whose existence is 
guaranteed by Proposition \ref{prop_adjunction}, so that 
$K_{\cal F}|_{S^{\nu}} = K_{\cal G}+\Delta$. We have
\[
\nu^*H_R\cdot (K_{\cal G}+\Delta)=\nu^*H_R\cdot \nu^*K_{\cal F}<0
\]
and so by applying bend and break  (e.g. see \cite[Corollary 2.28]{spicer20}), we may produce through any point $x \in S^\nu$ a rational curve $\Sigma$
with $\nu^*H_R\cdot \Sigma = 0$, a contradiction of the fact that $\loc R$ is one dimensional.
\end{proof}

\begin{lemma}
\label{l_alg_flip_cont}
Set up as above. 

Then the contraction of the flipping locus exists in the category of algebraic spaces.
\end{lemma}
\begin{proof}
By Lemma \ref{lem_null=loc}, it follows that $\Null H_R=\loc R$. 
Thus, Proposition \ref{p_existence_of_alg_contractions} implies the claim.
\end{proof}

\begin{remark}
\label{q_proj_flip_cont}
We remark that Lemma \ref{l_alg_flip_cont} holds equally well in the case where we only assume that $X$ is quasi-projective, $c\colon X \to Y$ is a contraction between quasi-projective varieties (or algebraic spaces)
and $R \subset \overline{NE}(X/Y)$ is a $K_{\mathcal F}$-negative extremal ray such that $\loc R$ is a curve $C$.
Indeed, to produce the contraction, we are free to replace $Y$ by an \'etale neighbourhood of $c(C)$ and so may assume that $Y$ is affine.  Further replacing $X$ and $Y$ by projective compactifications we may then apply Lemma \ref{l_alg_flip_cont} to produce the contraction.
\end{remark}

\subsection{Foliation complements}
\label{s_complement}
We  work in the following set up.
Let $X$ be a normal threefold  and let $\mathcal F$ be a 
foliation of rank one on $X$ with simple singularities (cf. Definition \ref{d_simple}). 
In particular, by Lemma \ref{l_simple_sing_are_quot}, 
$X$ admits at worst cyclic quotient singularities. 
Assume that $X$ admits a flipping contraction $f\colon X\to Z$ of a 
$K_{\cal F}$-negative connected curve $C$, where $Z$ is an algebraic space. 
Theorem \ref{t_cone2} and Remark \ref{r_cone}
imply that any component of  $C$ is $\mathcal F$-invariant and is not contained in $\sing \cal F$. 
Since $\cal F$ admits simple singularities, it follows that $X$ is smooth at any generic point of $C$.

We first consider the case that  $C$ is a smooth irreducible curve, whilst the case of a singular flipping curve will be considered in Section \ref{s_singularflippinglocus}.
The goal of this subsection is to  prove the following:

\begin{proposition}
\label{l_descend_vector_field}
Set up as above. 

Then, after possibly replacing $X$ by an analytic neighbourhood of $C$,
there exists a  divisor $T$ intersecting $C$ in a single point $Q$ such that
\begin{enumerate}
\item $(\cal F, T)$ is log canonical; 
\item $\cal F$ is terminal at $Q$; and
\item $K_{\cal F}+T \sim_{f,\mathbb Q} 0$.
\end{enumerate}
\end{proposition}

\begin{lemma} \label{l_singularities}
Set up as above.

Then there exists exactly one closed point $P\in C$ such that $\mathcal F$ is not terminal at $P$. 
Moreover $C\cap ( \sing X\cup \sing \cal F )$ consists of at most two points. 
\end{lemma}

\begin{proof}
Proposition  \ref{prop_reeb_stability} implies that there exists a point $P\in C$ such that $\mathcal F$ is not terminal at $P$. Let $Q\in C\cap \sing X$ be a closed point and assume, by contradiction, that $\cal F$ is terminal and Gorenstein at $Q$. Then Lemma \ref{l_terminal} implies that $C\subset \sing X$ and, in particular, the singularities  of $\cal F$ are not simple, a contradiction. 

Thus, since by assumption we have that $K_{\cal F}\cdot C<0$, the result follows immediately by Proposition \ref{prop_adjunction}. 
\end{proof}

\begin{lemma}
\label{l_restricted_contraction}
Set up as above.  Let $H$ be an ample divisor. 

Then, after possibly replacing $X$ by an analytic neighbourhood of $C$, there exists a divisor $L$ such that for any sufficiently large positive integer $m$ 
the general element $D$ of the linear system 
\[
\{ \Sigma \in |L+mH|\mid C\subset \Sigma\}
\]
is such that $D$ has at most two singularities along $C$, $D$ is klt  and $(D,C)$ is log canonical with a unique zero-dimensional log canonical centre along $C$. 

Moreover, if $f\colon X\to Z$ is the flipping contraction and $S$ is the normalisation of $f(D)$, then the induced morphism $f\vert_D\colon D\rightarrow S$ is a contraction of relative Picard number one.

\end{lemma}
\begin{proof}
The first part of the Lemma is a direct consequence of Lemma \ref{lem_preimagecurve}, Lemma \ref{l_generalelement} and  Lemma \ref{l_singularities}.

We now prove that the induced morphism $f|_D\colon D\to f(D)$ is a contraction of relative Picard number one.
By classical adjunction, we may write $(K_D+C)|_C = K_C+\Theta$
where $\Theta\ge 0$ is a $\mathbb Q$-divisor on $C$ which is supported on $\sing D\cap C$ and such that $(C,\Theta)$ is log canonical. 

Since there exists a unique zero-dimensional log canonical centre for $(D,C)$ 
along $C$, it follows  that 
the support of $\Theta$ consists of at most two points, of which only one of coefficient one for $\Theta$. 
Thus $(K_D+C)\cdot C <0$.  Since $f\vert_D$ only contracts the curve $C$ 
we see that $f\vert_D$ is in fact a $(K_D+C)$-negative contraction and is therefore of relative Picard number one.
\end{proof}

\begin{proposition}
\label{p_picard_rank_one}
Set up as above.  

Then the flipping contraction $f\colon X\to Z$ is a contraction 
of relative Picard number one.
\end{proposition}
\begin{proof}
This follows from Lemma
\ref{l_restricted_contraction}
and Lemma \ref{lem_pic_inj}.
\end{proof}

\begin{lemma}\label{l_terminal-point}
Set up as above. Suppose that $Q \in C$ is a  point where $\cal F$ is terminal
and $X$ is singular.

Then, after possibly replacing $X$ by an analytic neighbourhood of $C$, there exists an effective 
divisor  $T$ containing $Q$
 such that 
\begin{enumerate}
\item $(\cal F, T)$ is log canonical;
\item $K_{\cal F}+T$ is Cartier at $Q$; and
\item $(K_{\cal F}+T)\cdot C = 0$.
\end{enumerate}
\end{lemma}
\begin{proof}
%By Artin approximation Theorem, it suffices to produce $T$ on an analytic neighborhood of $C$.
Since $C$ is a curve we see that producing a divisor $T$ as required is in fact
an analytically local problem about $Q$. Thus, by Lemma \ref{l_terminal} and since $\cal F$ admits simple singularities, 
we may assume that there exists a cyclic quasi-\'etale morphism  $q\colon V \rightarrow X$
of order $m$, where $V\subset \bb C^3$ is an analytic open neighbourhood of the origin $0\in \mathbb C^3$, $q(0)=Q$  and 
the foliation  $\cal F'\coloneqq q^{-1}\cal F$ is induced by the  $\mathbb Z/m\bb Z$-equivariant  morphism 
\[
(x,y,z)\in V\mapsto (x,y)\in  \mathbb C^2.\]
By diagonalising this action we may freely assume that $\bb Z/m \bb Z$
acts by $(x, y,z) \mapsto (\zeta^a x, \zeta^b y, \zeta z)$ where $\zeta$
is a primitive $m$-th root of unity and $a,b$ are positive integers.
Note that $q^{-1}(C) = \{x = y = 0\}$.
Let $T' = \{z = 0\} \subset \bb C^3$ and let $T = q(T')$. We claim that $T$ satisfies all 
our desired properties.

First, $(\cal F', T')$ is clearly log canonical, and so it follows that $(\cal F, T)$
is log canonical by Lemma \ref{l_quasi-etale}.

Next, $T_{\cal F'}(-T')$ is generated by the vector field $z\frac{\partial}{\partial z}$ near $Q$
which is invariant under the  $\bb Z/m\bb Z$-action and therefore descends to a generating
section of $T_{\cal F}(-T)$. Thus, $K_{\cal F}+T$ is Cartier near $Q$.

Finally, by Lemma \ref{l_singularities} and Proposition \ref{prop_adjunction}, we have $K_{\cal F}\cdot C = -\frac{1}{m}$.  We claim that $T\cdot C = \frac{1}{m}$, 
from which our claim follows. Indeed, note that $T\cap C=\{Q\}$ and that $mT$ is Cartier at $Q$. 
Let $C'\coloneqq q^{-1}(C)$. 
Since $q|_{C'}\colon C' \rightarrow C$ is ramified to order $m$ at $Q$  and since $T'$ meets $C'$ transversally at one point, our claim follows. \end{proof}

\begin{proof}[Proof of Proposition \ref{l_descend_vector_field}] By Lemma \ref{l_singularities}, we have that if $\Sigma\coloneqq \sing X\cup \sing \cal F$, then $C\cap \Sigma$ consists of at most two points and it contains exactly one point at which $\cal F$ is not terminal. 
If $C\cap \Sigma$ contains two points, then by Lemma \ref{l_terminal-point} after possibly shrinking $X$ to an analytic neighbourhood of $C$, we may find a divisor $T$ 
 such that $(K_{\cal F}+T)|_C\sim_{\mathbb Q}0$ and $(\cal F,T)$ is log canonical. If $C\cap \Sigma$ consists of only one point then
 Proposition \ref{prop_adjunction} implies that  $K_{\cal F}\cdot C=-1$ and it follows immediately that there exists a divisor $T$, passing through a general point of $C$ and  satisfying the same properties as in the previous case.  
Thus, Proposition \ref{p_picard_rank_one} implies our claim.
\end{proof}

\subsection{Producing invariant divisors}
We work in  the same set up as in Section \ref{s_complement}. 
By Lemma  \ref{l_singularities}, there exists a unique closed point 
$P\in C$ at which $\cal F$ is not terminal. 
The goal of this section is to provide a precise description of the neighbourhood of a flipping curve, and use this precise description to produce a large number of $\cal F$-invariant divisors containing $C$.

\begin{proposition}
\label{prop_first_integral}
Set up as above. 

Then, in an analytic neighbourhood $U$ of $C$ there exists a projective variety $W$ 
and  a meromoprhic 
map $F\colon U \dashrightarrow W$ which is holomorphic on $U \setminus C$ 
such that $\mathcal F$ is induced 
by $F$.

Moreover, 
\begin{enumerate}
\item 
$X$ is smooth at $P$;
\item the semi-simple part of a vector field 
defining $\mathcal F$ near $P$ has eigenvalues $1, -a, -b$ where $a, b \in \mathbb Q_{>0}$; and 
\item there exists a $\mathcal F$-invariant $\mathbb Q$-divisor $D\ge 0$
such that $(U, D)$ is log canonical and $C$ is a log canonical centre of $(U, D)$.
\end{enumerate}
\end{proposition}
\begin{proof}
Let $T$ be the divisor whose existence is guaranteed by Proposition \ref{l_descend_vector_field}. 
Let $\cal G$ be the induced foliation on $Z$ and let  $D = f_*T$. 
Since $K_{\cal F}+T=f^*(K_{\cal G}+D)$, 
 we have that $(\cal G, D)$ is log canonical. After 
 replacing $Z$ by a quasi-\'etale cover of $Z$, we may assume that $K_{\cal G}+D$ 
is Cartier and $\cal G(-D)$ is generated by a vector
field ${\partial}$.  Consider an embedding $\iota\colon Z \hookrightarrow \bb C^m$ and a lift $\widetilde{\partial}$
of $\partial$ to a vector field on $\bb C^m$. 
%Let $m$ denote the dimension of the Zariski tangent space to $Z$ at $R := f(C)$.

Propostion \ref{p_semisimple} implies that, up to a formal change of coordinates and rescaling, 
$\widetilde{\partial}$ is a semi-simple vector field 
and
$\widetilde{\partial} = \sum_{i = 1}^{m'} \lambda_ix_i\partial_{x_i}$
where $m'\leq m$ and $\lambda_1, \dots, \lambda_{m'}$ are positive integers.
We may apply a theorem of Poincar\'e
(see \cite[Remarques historiques 3.3]{martinet81}) to see that we may in fact take this change 
of coordinates to be holomorphic.
We take $U$ to be the pre-image under $f$ of the neighbourhood of $f(C)$ where this coordinate change is well defined.

Let $\mathcal H$ denote the foliation induced by $\widetilde{\partial}$.
Let $b\colon \overline{\mathbb C^m} \rightarrow \mathbb C^m$ be the weighted blow up 
in $x_1, \dots, x_{m'}$ with weights $\lambda_1, \dots, \lambda_{m'}$.
It is easy to check that $b^{-1}\mathcal H$ admits a holomorphic first integral
$\Phi\colon \overline{\mathbb C^m} \rightarrow \mathbb P(\lambda_1, \dots, \lambda_{m'}) \times \mathbb C^{m-m'}$.  This induces 
a meromoprhic map $F\colon X \dashrightarrow \mathbb P(\lambda_1, \dots, \lambda_{m'}) \times \mathbb C^{m-m'}$ which is a
meromorphic first integral of $\mathcal F$.

Since $\mathcal G$ has canonical singularities away from $R:=f(C)$, we see that 
$\Phi\vert_Z$ is holomorphic on $Z \setminus R$, and hence $F$ is holomorphic
on $X \setminus C$. 

\medskip

We now verify our three remaining claims.

We first show (1).
Assume for sake of contradiction that  $X$ is not smooth at $P$. 
Since $\cal F$ admits simple singularities, 
there exists  an analytic open neighbourhood $V$ of $P$
such that the restriction of $\cal F$ on $V$ is as in Example \ref{ex_unresolvable}. 
In particular, $K_{\cal F}$ is not Cartier at $P$. On the other hand, 
we have that $K_{\cal F}+T$ is Cartier and,  Proposition \ref{l_descend_vector_field} implies that
$T$ intersect $C$ in a single point $Q$ such that $\cal F$ is terminal at $Q$. 
In particular, $Q\neq P$ and therefore $K_{\cal F}$ is Cartier at $P$, a contradiction. 
Thus, $X$ is smooth at $P$.

We now show (2). We observe that 
the conditions of Lemma \ref{lem_eigenvalue_computation}
are satisfied by $C$ and $f^*x_1, \dots, f^*x_m$, and so we may apply the Lemma to conclude.

Finally we verify (3).
Let $\overline{Z}$ be the strict transform of $Z \subset \mathbb C^m$ under $b$, let $\overline{X}$ be the normalisation of the component of $X\times_Z\overline{Z}$ which dominates $Z$
 and let $ \overline F\colon \overline{X} \to \mathbb P(\lambda_1, \dots, \lambda_{m'}) \times \mathbb C^{m-m'}$
be the composition of the projection $\pi\colon \overline{X} \to \overline{Z}$ with restriction of $\Phi$ to $\overline{Z}$.
Notice that we have a birational contraction $p\colon \overline{X} \to X$ which defines an isomorphism $\overline{X} \setminus \exc \pi \to X \setminus C$.  Moreover, $ \overline F$ yields a holomorphic first integral of $p^{-1}\mathcal F$.

Let $A$ be an ample divisor on $\mathbb P(\lambda_1, \dots, \lambda_{m'}) \times \mathbb C^{m-m'}$ and let $H \in |kA|$ be a general element, where $k \gg 0$.
Since $p^{-1}\mathcal F$ has simple 
singularities on $\overline{X} \setminus \exc p$, we deduce that $(\overline{X} \setminus \exc p,  \overline F^*H|_{\overline{X}\setminus \exc p})$
is a simple normal crossings pair.  In particular, 
$(X \setminus C, p_* \overline F^*H|_{X \setminus C})$ is log canonical.
Since $p(\exc p) = C$, by taking $k$ to be sufficiently large, the multiplicity of the divisor $p_* \overline F^*H$ along $C$ can be made arbitrarily large and so $(X, p_* \overline F^*H)$ will not be log canonical 
at the generic point of $C$.  

Let $\lambda$ be the log canonical threshold of $X$ with respect to $p_* \overline F^*H$
and set $D := \lambda p_* \overline F^*H$.  
Then $C$ is a log canonical centre of $(X, D)$ and 
$(X, D)$ is log canonical away from finitely many closed points of $X$. 
Theorem \ref{t_isolated_sing_bound} then applies to show that $(X, D)$ is log canonical, 
and we may conclude.
\end{proof}

\subsection{Singular flipping locus}\label{s_singularflippinglocus}

We now show that if $X$ is a normal threefold and $\cal F$ is a foliation of rank one on $X$ with simple singularities and which admits a flipping contraction $f\colon X\to Z$ of a $K_{\cal F}$-negative irreducible curve $C$ then $C$ is a smooth curve. Our method was inspired by \cite[II.i]{mcq04}. 
We begin with the following:

\begin{lemma}
\label{lem_no_flipping_tangencies}
Let $\partial$ be a vector field defined over a neighbourhood of 
$0 \in \bb C^3$ and assume that, in suitable coordinates, we may write
\[
\partial = at\frac{\partial}{\partial t}-b x\frac{\partial}{\partial x}-c y\frac{\partial}{\partial y}
\]
where $a, b, c$ are positive integers.
Let $C = \{x = y = 0\}$
and $D$ be a  $\partial$-invariant prime divisor such that $D\cap C = \{0\}$. 

Then $D$ meets $C$ transversely.
\end{lemma}
\begin{proof}
We may write $D = \{f = 0\}$ where $f$ is $(a, -b, -c)$-weighted homogeneous of degree $d$, i.e. 
\[
f(t,x, y) = \sum_{ia-bj-ck = d} a_{ijk}t^ix^jy^k
\]
for some $a_{ijk}\in \mathbb C$. 
 Since
$D$ does not contain $C$ we see that $f$ is not an element of the ideal $(x, y)$, 
which implies that $a_{i00}$ is non-zero for some $i>0$. 
In particular, $d$ is a positive integer and, therefore, 
$a_{0jk} = 0$ for all $j, k\ge 0$.  Thus, $D = \{t = 0\}$ and our result follows.
\end{proof}

\begin{proposition}\label{p_notsingular}
Let  $X$ be a normal threefold and let $\cal F$ be a foliation of rank one on $X$ with simple singularities and which admits a flipping contraction $f\colon X\to Z$ of a $K_{\cal F}$-negative irreducible curve $C$. 

Then $C$ is a smooth curve. 
\end{proposition}

\begin{proof}
Suppose by contradiction that $C$ is not smooth. As in the proof of Lemma \ref{l_singularities}, 
Proposition \ref{prop_adjunction} implies that $C$ admits a unique cusp at a point $P\in C\cap \sing \cal F$.  We first prove the following:

\begin{claim} 
There exists a birational morphism $p\colon X' \rightarrow X$ such that if $\cal F' \coloneqq p^{-1}\cal F$ and $C'$ is the strict transform of $C$ in $X'$ then 
\begin{enumerate}
\item $C'$ is smooth;

\item there is a $p$-exceptional prime divisor $E$ in $X'$ which is $\cal F'$-invariant and is tangent to $C'$;

\item $K_{\cal F'} = p^*K_{\cal F}$; and

\item $C'$ spans a $K_{\cal F'}$-negative extremal ray $R'$.
\end{enumerate}

\end{claim}
\begin{proof}[Proof of the Claim.]
Lemma \ref{lem_preimagecurve} implies that $X$ is smooth at $P$.
We may find a sequence of blow ups 
\[p\colon X'=X_n\stackrel{p_{n}}{\longrightarrow}X_{n-1}\longrightarrow\dots \longrightarrow X_1\stackrel {p_1}\longrightarrow X\] 
in $\cal F$-invariant closed points which resolve the cusp of $C$ at $P$.
Let $E$ be the  $p_n$-exceptional divisor in $X'$ and let $C'$ be the strict transform 
of $C$ in $X'$. We may assume that $p_n(C')$ is singular, which implies that $E$ is tangent to $C'$. 
Let $\cal F'=p^{-1}\cal F$.
Lemma \ref{l_nondicritical} implies that $E$ is $\cal F'$-invariant. 
 By \cite[Lemma I.1.3]{bm16}, we have that $K_{\cal F'} = p^*K_{\cal F}$.

To prove (4), let $G$ be a $p$-exceptional divisor so that $-G$ is $p$-ample and let $H_R$ be
the supporting hyperlane of the ray $R$ spanned by $C$.
Then for $\delta>0$ sufficiently small we may find an ample divisor $A$ on $X'$
so that $p^*H_R-\delta G+A$ is a big and nef divisor which is only zero on the strict transform of
curves in $\Null H_R$. Thus, $C'$ spans a $K_{\cal F'}$-negative extremal ray, as claimed. 
\end{proof}

We now proceed with the proof of the Proposition. 
We may apply Lemma \ref{l_alg_flip_cont} (cf. Remark \ref{q_proj_flip_cont}) to see that there exists a flipping contraction $f'\colon X'\to Z'$  in the category of algebraic spaces associated to $R'$.
Let $P' = C' \cap \sing \cal F'$
and let $\partial'$ be a local generator of $\cal F'$ near $P$.
By Proposition \ref{prop_first_integral}.(2), 
after a suitable renormalisation,  the semi-simple part of $\partial'$
has eigenvalues $(a, -b, -c)$ where $a, b, c$ are all positive integers.
Thus,  Lemma \ref{lem_no_flipping_tangencies} implies that $E$ is 
transverse to $C'$, a contradiction.
\end{proof}

We now show that each connected component of  the flipping locus is irreducible. The same result may be found in \cite{mcq04}.

\begin{lemma}
\label{lem_meets_flip}
Let $X$ be a normal threefold
and let $\cal F$ be a rank one foliation with simple  singularities.
Let $c\colon X\to Y$ be a projective morphism in the category of algebraic spaces
and let $C_1$ and $C_2$ be two distinct irreducible curves in $X$ such that  $C_1 \cap C_2 \neq \emptyset$.
Assume that $R_1=\mathbb R_+[C_1]$ and $R_2=\mathbb R_+[C_2]$ are distinct $K_{\mathcal F}$-negative
extremal rays of $\overline{NE}(X/Y)$. 
Suppose furthermore that $\loc(R_1)=C_1$ and that 
the flipping contraction and flip associated to $R_1$ exist.

Then for a general $x\in X$, there exists a $\mathcal F$-invariant curve $\Sigma_x$ in $X$ passing through $x$ and  rational numbers $a,b\ge 0$ such that $[aC_1+bC_2]=[\Sigma_x]$ in 
$\overline{NE}(X/Y)$.
\end{lemma}
\begin{proof}
Consider the flip $\phi\colon X \dashrightarrow X'$ of $C_1$ and let $C'_2$
be the strict transform of $C_2$ in $X'$.  It follows from the negativity lemma (cf. Lemma \ref{l_negativity})
that if $\cal F':=\phi_*\cal F$ then $\cal F'$ is terminal at all, not necessarily closed, points of $C'_2$. 
By Proposition  \ref{prop_reeb_stability}, we may assume  that there exists a point $P\in C_2$ such that $\mathcal F$ is not terminal at $P$. As in the proof of Lemma \ref{l_singularities}, it follows that $C_2\cap ( \sing X\cup \sing \cal F )$ consists of at most two points. Thus, there are at most two  terminal 
non-Gorenstein singularities along $C'_2$ and so we may apply foliation adjunction (cf. 
Proposition \ref{prop_adjunction}) to deduce that  $K_{\cal F'}\cdot C'_2 <0$.
Therefore, Proposition \ref{prop_reeb_stability} implies 
that $C'_2$ moves in a family of $\cal F'$-invariant curves. Thus, the claim follow. 
\end{proof}

%% file: section6.tex
\section{Threefold contractions and flips}

\subsection{Divisorial contractions}

\begin{lemma}
\label{lem_lc_contraction} 
Let $X$ be a $\mathbb Q$-factorial klt projective threefold and let $\cal F$ be a rank one foliation 
on $X$ with  canonical singularities. Let $R$ be a $K_{\cal F}$-negative extremal ray such that 
 $D\coloneqq \loc R $ has  dimension two. 
 
 Then
 \begin{enumerate}
 \item $D$ is $\mathcal F$-invariant; and
 \item  
 if $\Gamma \geq 0$ if a $\mathbb Q$-divisor on $X$ with $\cal F$-invariant support and such that $(X, \Gamma+D)$ is log canonical, then the divisorial contraction $c_R\colon X\to Y$ associated to $R$ exists in the category of projective varieties.
 
 \end{enumerate}

\end{lemma}
\begin{proof}
Note that $D$ is an irreducible divisor. Let $\nu\colon D^\nu \rightarrow D$ be the normalisation.
 and suppose for the sake of contradiction that $D$
is not $\mathcal F$-invariant.  

Let $H_R$ be the supporting hyperplane to $R$.  By Lemma \ref{lem_nakamaye} 
we have for any ample divisor $A$ and $\epsilon>0$ sufficiently small that $\mathbb B(H_R-\epsilon A) = D$.
In particular, if $m>0$ is sufficiently divisible we may write $m(H_R-\epsilon A) = kD+G$ where $k>0$
and $G$ is movable.  In particular, it follows that $\nu^*D \sim_{\mathbb Q} \frac{1}{k}(m(H_R-\epsilon A)-G)$
is not pseudo-effective.  From this we conclude that $\nu^*(K_{\mathcal F}+D)$ is not pseudo-effective.
On the other hand, by foliation adjunction, \cite[Proposition-Definition 3.7]{CS23b}
$\nu^*(K_{\mathcal F}+D) \sim_{\mathbb Q} \Delta \ge 0$, a contradiction.

\medskip
We will now show that the contraction exists supposing 
that  $\Gamma \geq 0$ if a $\mathbb Q$-divisor on $X$ with $\cal F$-invariant support and such that $(X, \Gamma+D)$ is log canonical.
We will prove that $R$ is $(K_X+\Gamma+D)$-negative. 
Let $\cal G$ be the  foliation on $D^{\nu}$ and $\Delta$ be the $\mathbb Q$-divisor, whose existence is guaranteed by  Proposition \ref{prop_adj_comp}
and let $\Theta\ge 0$ be the $\mathbb Q$-divisor on $D^{\nu}$ such that
\[
(K_X+\Gamma+D)|_{D^\nu} = K_{S^\nu}+\Theta.
\]
Since $D$ is covered by curves $\xi$ such that $(K_{\cal G}+\Delta )\cdot \xi<0$, 
by a similar argument as in the proof of  Lemma \ref{lem_KX+inv_antinef}, it follows that $\mathcal G$ is algebraically integrable. 
Proposition \ref{prop_adj_comp} also implies that for any curve $C\subset D^{\nu}$ which is not 
$\cal G$-invariant, we have that $\mu_C\Delta \geq \mu_C\Theta$.  
Since $(D^\nu, \Theta)$ is log canonical and since $\cal G$ is algebraically integrable,  Lemma \ref{lem_KX+inv_antinef2} implies 
 that all the $(K_{\cal G}+\Delta)$-negative curves in $D^\nu$ which are $\cal G$-invariant  
are in fact $(K_{D^\nu}+\Theta)$-negative. 
Thus, $R$ is $(K_X+D)$-negative and, therefore, the divisorial contraction associated to $R$ exists \cite[Theorem 5.6]{ambro01}.\end{proof}

\begin{theorem}\label{t_div_contraction}
Let $X$ be a projective $\mathbb Q$-factorial klt threefold and let $\cal F$ be a rank one foliation on $X$ with  canonical singularities. Let $R$ be a $K_{\cal F}$-negative extremal ray such that 
 $D\coloneqq \loc R $ has  dimension two. Let $\Gamma \geq 0$ be a $\mathbb Q$-divisor on $X$ with $\cal F$-invariant support, and such that $D$ is not contained in the support of $\Gamma$ and $(X, \Gamma)$ is log canonical.

Then the divisorial contraction  associated to  $R$ exists.  In particular,  there exists a projective birational morphism $c_R\colon X \rightarrow Y$, whose exceptional divisor coincides with $D$ and such that, if $\cal F'$ is the foliation induced on $Y$ then 

\begin{enumerate}
\item $Y$ is projective; 
\item $\rho(X/Y) = 1$; 
\item $\cal F'$ has canonical singularities and it is terminal at every point of  $c(\exc c )$; and

\item $(Y, (c_R)_*\Gamma)$ is log canonical.
\end{enumerate}
\end{theorem}
\begin{proof}
If $(X, \Gamma+D)$ is log canonical we may apply Lemma \ref{lem_lc_contraction} to produce our desired contraction.

So assume that $(X, \Gamma+D)$ is not log canonical.
Let $\lambda$ denote the log canonical threshold  of $X$ with respect to $D$. Then $\lambda<1$ and Theorem \ref{t_isolated_sing_bound} implies
that $(X,\Gamma+\lambda D)$  admits a one-dimensional log canonical centre $C\subset X$. 
 Proposition \ref{p_boundingsing} implies that $C$ is not contained in $\sing \cal F$. 
Let $\nu\colon D^{\nu}\to D$ be the normalisation of $D$.
By Proposition \ref{prop_adjunction}, there exists a foliated pair $(\cal G,\Delta)$ on $D^\nu$ such that 
\[K_{\cal F}|_{D^\nu}=K_{\cal G}+\Delta. 
\]
	\begin{claim}
		$C$ is $\cal F$-invariant.
	\end{claim}
	\begin{proof}
		By \cite[Lemma 4.2]{CS23b} to check invariance we may freely replace $X$ by the index one cover
		associated to $K_{\mathcal F}$ in a neighbourhood of a general point of $C$.
		Since $(X, D)$ is not log canonical it follows that $C \subset \sing X \cup \sing D$, and so 
		by \cite[Theorem 5]{MR0212027} we conclude that $C$ is $\cal F$-invariant.
	\end{proof}
	Since $C$ is not contained in $\sing \cal F$ and $\nu^{-1}(C)$ is not contained in the singular locus of $D^{\nu}$, it follows that $\nu^{-1}(C)$ is $\cal G$-invariant.
Since $D=\loc R$, it follows that  $[C]\in R$  and, in particular, $K_{\cal F}\cdot C<0$.  
 Theorem \ref{thm_subadj} implies that $(K_X+\Gamma+\lambda D)\cdot C<0$ and so $R$
is $(K_X+\Gamma+\lambda D)$-negative.  Thus, we can realise the $K_{\cal F}$-contraction as a $(K_X+\Gamma+\lambda D)$-negative contraction. In particular, (1) and (2) hold. 
Lemma \ref{l_negativity} implies (3). 
%
%We now prove (4). Note that, using the same argument as above, we may take the contraction  to be associated to the log pair $(X, \Gamma+\lambda D)$
%where $\lambda \geq 0$ is the log canonical threshold of $(X, \Gamma)$ with respect to $D$. Thus, 
The  negativity lemma (cf. \cite[Lemma 3.38]{KM98})  implies (4). 
\end{proof}

\subsection{Flips}

\begin{lemma}
\label{l_flips_exist} 
Let $X$ be a normal  threefold and 
let $\cal F$ be a rank one foliation on $X$ with  simple singularities. 
Let $c\colon X\to Y$ be a projective morphism in the category of algebraic spaces and 
let $R$ be a $K_{\cal F}$-negative extremal ray of  $\overline{NE}(X/Y)$ such that 
 $\loc R $ has dimension one and $c\colon X\to Y$ is the associated flipping contraction. 
Let $H_R$ be a supporting hyperplane to $R$ for $\overline{NE}(X/Y)$.

 Then each connected component of $\exc c$ is irreducible,  the flip  associated to $R$ exists and  $H_R$ 
descends to a $\bb Q$-Cartier divisor $M$ on $Y$.
\end{lemma}

\begin{proof}
Lemma \ref{l_simple_sing_are_quot} implies that $X$ has quotient singularities. 
In particular, $X$ is klt and $\mathbb Q$-factorial. 
The problem of  descending $H_R$ and of 
constructing the flip is \'etale local on the base. Thus,  we may freely replace
$Y$ by an \'etale neighbourhood of a point in $c(\exc c)$. 

By shrinking about a Zariski neighbourhood of $c(\exc c)$ we may freely assume that $\exc c$ is connected.
We will show that $\exc c$ is in fact irreducible and that the flip exists.  Let $C_1, \dots, C_r$ be the irreducible components 
of $\exc c$.  

We first claim that after replacing $Y$ by an \'etale neighbourhood of $c(\exc c)$,
we may assume that $C_1, \dots, C_r$ span distinct extremal rays in $\overline{NE}(X/Y)$.
Indeed, let $\widehat{X}$ denote the formal completion of $X$ along $\exc c$ and let $\hat{c}$ denote the restricted map. 
Then, for any $i=1,\dots,r$,  we may
find a formal $\bb Q$-Cartier divisors $D_i \subset \widehat{X}$ such that $D_i\cdot C_j = \delta_{ij}$ for any $j=1,\dots,r$, where $\delta_{ij}$ is the Kronecker delta.  By the approximation theorems (cf. Section \ref{s_approx}), 
after replacing $Y$ by an \'etale neighbourhood of $c(\exc c)$, for any $j=1,\dots,r$,
we may find a divisor  $\tilde D_j$
%on an \'etale neighbourhood of $c(\exc c) \in Y$ 
which approximate $\hat{c}_*D_{j}$. Thus, our claim follows.
%This will be our desired \'etale neighbourhood.

Let $R_1=\mathbb R_+[C_1]$. 
By Lemma \ref{l_alg_flip_cont} (cf. Remark \ref{q_proj_flip_cont}) the contraction $f\colon X \to Z$ over $Y$ associated to $R_1$ exists.
We will show that the flip of $R_1$ exists.  
Let $D\ge 0$ be a $\cal F$-invariant $\mathbb Q$-divisor in an analytic neighbourhood of $C_1$ such that $(X,D)$ is log canonical around $C_1$ and $C_1$ is a log canonical centre of $(X,D)$ and whose existence is guaranteed by Proposition \ref{prop_first_integral}(3).
Theorem \ref{thm_subadj} implies that $(K_X+D)\cdot C <0$.

Fix $n\geq 0$ and let $X_n$ denote the $n$-th infinitessimal neighborhood of $C_1$ in $X$.
By our approximation results (cf.  Section \ref{s_approx}), after possibly replacing  
$Z$ by an \'etale neighborhood of $f(C_1)$, we may 
find a divisor $\tilde{D}$ 
such that $\tilde{D}\vert_{X_n} = D\vert_{X_n}$.
By Lemma \ref{l_formalklt}, it follows that taking $n$ to be sufficiently large, the pair $(X, \tilde D)$ is log canonical
and  
\[
(K_X+\tilde D)\cdot C<0.
\]
In particular,  the $K_{\cal F}$-flipping contraction (resp. flip) can be realised as a 
$(K_X+\tilde D)$-flipping contraction (resp. flip) and the basepoint free theorem implies that 
$H_R$ descends to a $\mathbb Q$-Cartier divisor on $Z$.

We may now apply Lemma \ref{lem_meets_flip} to see that in fact $\exc c$ is irreducible, hence
$Z = Y$ and the flip of $R_1$ is in fact the flip of $R$.
\iffalse
By   Proposition \ref{p_notsingular} and Lemma \ref{lem_meets_flip}, 
we may assume that $R=\mathbb R_+[C]$, where $C$ is a smooth irreducible curve. 
Let $D\ge 0$ be a $\cal F$-invariant $\mathbb Q$-divisor in an analytic neighbourhood of $C$ such that $(X,D)$ is log canonical around $C$ and $C$ is a log canonical centre of $(X,D)$ and whose existence is guaranteed by Proposition \ref{prop_first_integral}.(3).
Theorem \ref{thm_subadj} implies that $(K_X+D)\cdot C <0$.

Fix $n\geq 0$ and let $X_n$ denote the $n$-th infinitessimal neighborhood of $C$ in $X$.
By our approximation results (cf.  Section \ref{s_approx}), after possibly replacing  
$Y$ by a smaller \'etale neighborhood, we may 
find a divisor $\tilde{D}$ 
such that $\tilde{D}\vert_{X_n} = D\vert_{X_n}$.
By Lemma \ref{l_formalklt}, it follows that taking $n$ to be sufficiently large, the pair $(X, \tilde D)$ is log canonical
and  
\[
(K_X+\tilde D)\cdot C<0.
\]
In particular,  the $K_{\cal F}$-flipping contraction (resp. flip) can be realised as a 
$(K_X+\tilde D)$-flipping contraction (resp. flip) and the basepoint free theorem implies that 
$H_R$ descends to a $\mathbb Q$-Cartier divisor on $Y$, as claimed. 
\fi
\end{proof}

\begin{theorem}\label{thm_flips_exist}Let $X$ be a normal projective threefold and 
let $\cal F$ be a rank one foliation on $X$ with  simple singularities. 
Let $R$ be a $K_{\cal F}$-negative extremal ray such that 
 $ \loc R $ has  dimension one.

Then the flipping contraction  $c_R\colon X\to Y$ associated to $R$  
exists in the category of projective varieties.  
Moreover, the flip $\phi\colon X \dashrightarrow X^+$ associated to $R$ exists 
and if $\cal F^+$ is the foliation induced on $X^+$ then 

\begin{enumerate}
\item $X^+$ is projective and has quotient singularities; 

\item $\rho(X/Y) = \rho(X^+/Y)=1$;

\item $\cal F^+$ has simple singularities and 
$\cal F^+$ is terminal at every point of $\exc \phi^{-1}$; and

\item if $\Gamma \geq 0$ is a $\mathbb Q$-divisor on $X$ with $\cal F$-invariant support such that $(X, \Gamma)$ is log canonical, then
$(X^+,\phi_* \Gamma)$ is log canonical.
\end{enumerate}

\end{theorem}
\begin{proof}
Lemma \ref{l_simple_sing_are_quot} implies that $X$ has quotient singularities. 
In particular, $X$ is klt and $\mathbb Q$-factorial. 
Let $c_R\colon X\to Y$ be the flipping contraction associated to $R$ 
in the category of algebraic spaces and whose existence is guaranteed by Lemma \ref{l_alg_flip_cont}.
Let $H_R$ be a supporting hyperplane to $R$ for $\overline{NE}(X)$.
By Lemma \ref{l_flips_exist}, each connected component of $\exc c_R$ is irreducible, 
$H_R$ 
descends to a $\bb Q$-Cartier divisor $M$ on $Y$ and the flip $\phi\colon X\dashrightarrow X^+$ associated to $R$ exists. In particular, 
$M^{\dim Z}\cdot Z >0$ for all positive dimensional $Z \subset Y$ and so $M$ is ample
by the Nakai-Moishezon criterion and, in particular, $Y$ is projective.

%To prove the projectivity of $Y$ it suffices to show that $H_R$ 
%descends to a $\bb Q$-Cartier divisor $M$ on $Y$.  Indeed, 
%in this case $M^{\dim Z}\cdot Z >0$ for all positive dimensional $Z \subset Y$ and so $M$ is ample
%by the Nakai-Moishezon criterion and, in particular, $Y$ is projective.
%
%
%The problem of descending $H_R$ and of constructing the flip is analytic local on the base. Thus,  we may freely replace
%$Y$ by an analytic neighbourhood of a point in $c_R(\exc c_R)$. 
%
%Let $D$ be a $\cal F$-invariant divisor in an analytic neighbourhood of $C$ such that $(X,D)$ is log canonical around $C$ and $C$ is a log canonical centre of $(X,D)$ and whose existence is guaranteed by Corollary \ref{cor_flipping_lc_centre}.
%Theorem \ref{thm_subadj} implies that $(K_X+D)\cdot C <0$.
%
%
%Fix $n\geq 0$ and let $X_n$ denote the $n$-th infinitessimal neighborhood of $C$ in $X$.
%By our approximation results (cf.  Section \ref{s_approx}), after possibly replacing  
%$Y$ by a smaller \'etale neighborhood, we may 
%find a divisor $\tilde{D}$ 
%such that $\tilde{D}\vert_{X_n} = D\vert_{X_n}$.
%By Lemma \ref{l_formalklt}, it follows that taking $n$ to be sufficiently large, the pair $(X, \tilde D)$ is log canonical
%and  
%\[
%(K_X+\tilde D)\cdot C<0.
%\]
%In particular,  the $K_{\cal F}$-flipping contraction (resp. flip) can be realised as a 
%$(K_X+\tilde D)$-flipping contraction (resp. flip) and the basepoint free theorem implies that 
%$H_R$ descends to a $\mathbb Q$-Cartier divisor on $Y$.

%From this we deduce that $Y$ and $X^+$ are projective and
Thus, also $X^+$ is projective and 
 $\rho(X/Y)=\rho(X^+/Y)=1$. By Proposition \ref{prop_admiss_sing_preserved}, it follows that $\cal F^+$ has simple singularities, and 
Lemma \ref{l_simple_sing_are_quot} implies that $X^+$ has quotient singularities. Thus, (1) and (2) follow. 
Lemma \ref{l_negativity} implies (3).

We now prove (4). Let $\Gamma$ be an $\cal F$-invariant divisor such that 
$(X, \Gamma)$ 
is log canonical.  
As in the proof of Proposition \ref{prop_first_integral}(3), up to replacing $X$ by an analytic neighbourhood of a connected component $C$ of $\exc c_R$,
we may find a $\mathbb Q$-divisor $D\ge 0$ whose support is $\cal F$-invariant and such that $(X, \Gamma+D)$ is not log canonical and 
%and, if $\widehat X$ denotes the formal completion of $X$ at the generic point of $C$ then 
 $C$ is the only non-log canonical centre of $(X,\Gamma+D)$ of positive dimension. Thus,  if $\lambda$ is the log canonical threshold of $( X,\Gamma)$ with respect to $D$ along $C$ then by Theorem \ref{t_isolated_sing_bound} we have that  $( X, \Gamma+\lambda D)$ is log canonical and by Theorem \ref{thm_subadj}, we have that  
 $-(K_{ X}+\Gamma+\lambda D)$ is ample over $Y$.  It follows by the negativity lemma 
 (cf. \cite[Lemma 3.38]{KM98}) that $(X^+, \phi_*(\Gamma+\lambda D))$ is log canonical
and, therefore, $(X^+, \phi_*\Gamma)$ is log canonical. Thus, (4) follows. 
\end{proof}

%% file: section7.tex
\section{Termination of flips}
The goal of this section is to prove the following: 

\begin{theorem}[Termination of flips]
\label{thm_term}
Let $X$ be a normal variety and let $\cal F$
be a rank one foliation on $X$ with canonical singularities.

Then any sequence of $K_{\cal F}$-flips terminates.
\end{theorem}

We begin with the following
\begin{lemma}
\label{l_sing_trans}
Let $X$ be a normal  variety and let $\cal F$
be a rank one foliation on $X$ with canonical singularities.
Let $\phi\colon X \dashrightarrow X^+$ be a $K_{\cal F}$-flip and let $Z^+\subset X^+$ be the flipped locus.

Then $Z^+\cap \sing \cal F^+ = \emptyset$.
\end{lemma}

Note that the corresponding statement for higher rank foliations, including the absolute case, is easily shown to be false.

\begin{proof}
Suppose not and let $P\in Z^+\cap \sing \cal F^+$ be a closed point. Then Lemma \ref{l_terminal}  implies that $\cal F^+$ is not terminal near $P$. 
Thus, there exists an exceptional divisor $E$ over $X$ centred at $P$ and such that $a(E,\cal F^+)=0$. The negativity Lemma (cf. Lemma \ref{l_negativity}) implies that $a(E,\cal F)<0$, a contradiction. 
\end{proof}

\begin{proposition}[Special termination]
\label{prop_special_term}
Let $X$ be a normal  variety and let $\cal F$
be a rank one foliation on $X$ with canonical singularities.
Let
\[X=X_0  \dashrightarrow X_1 \dashrightarrow X_2\dashrightarrow \dots \]
be a sequence of $K_{\cal F}$-flips and let $\cal F_i$ be the induced foliation on $X_i$.  

Then, after finitely many flips,
the flipping and flipped locus do not meet any log canonical
centres of $\cal F_i$ properly.
\end{proposition}

Note that, using the same notation as in Proposition \ref{prop_special_term}, since $\cal F_i$ is canonical, a log canonical centre for $\cal F_i$ is just a canonical centre. Moreover, by Lemma \ref{l_terminal}, if $P\in X$ is a zero-dimensional log canonical centre for $\cal F$ then $P\in \sing^+\cal F$.

\begin{proof} Let $\phi_i\colon X_i\dashrightarrow X_{i+1}$ denote the $K_{\cal F_i}$-flip and let  $S_i \coloneqq \sing \cal F_i$.
By Lemma \ref{l_sing_trans}, it follows that $\phi_i^{-1}$ is isomorphic around $S_{i+1}$. Therefore, the number of irreducible components of $S_i$ is not increasing as $i$ increases. 

Lemma \ref{l_sing_trans} also implies that if a connected component of $S_i$ intersects the flipping locus, then it is contained in the flipping locus and, therefore,  the number of connected components of $S_i$ decreases after such a flip. Thus, our claim follows. 
%By Lemma \ref{l_sing_trans}, it follows that $\phi_i$
%is either an isomorphism around $S_i$ 
%or it contracts some irreducible component of $S_i$.  Since the latter can
%happen only finitely many times, we may assume
%that, after finitely many steps,  each $\phi_i$ is an isomorphism around $S_i$. Let $\psi_i\colon S_i\to S_{i+1}$ denote the induced isomorphism. 
%
%We now show that the flipping locus is disjoint
%from $S_i$.  Let $T_i \subset S_i$ denote the intersection
%of the flipping locus with $S_i$ and assume it is not empty.  Observe
%that $\psi_i(T_i) \subset S_{i+1}$ is contained in the intersection 
%of the flipped locus with $S_{i+1}$, but this implies that the flipped locus meets $S_{i+1}$, a contradiction. 
\end{proof}

\begin{remark}
In fact, this argument shows that each flip contracts an entire
component of the singular locus of the foliation, i.e., 
if $Z \subset \sing \cal F$ meets the flipping locus
then in fact it is contained in the flipping locus.  
This also follows from the explicit description
of the flip given in \cite{mcq04}, 
but it is interesting to note that this can also be proven
by a simple discrepancy calculation.
\end{remark}

\begin{proof}[Proof of Theorem \ref{thm_term}]
By Lemma \ref{l_terminal} and  Proposition \ref{prop_reeb_stability},
 it follows  that if $C \subset X$ is a flipping
curve then $C$ must meet $\sing \cal F$ at some point and, in particular,
it meets a log canonical centre of $\cal F$.  Thus,
Proposition \ref{prop_special_term} implies the claim. \end{proof}

\section{Running the MMP}

\subsection{Running the MMP with simple singularities}

\begin{proposition}
\label{prop_admissible_mmp}
Let $X$ be a normal projective threefold and let $(\cal F,\Delta)$ be a rank one foliated pair on $X$ with log canonical singularities and such that $\cal F$ admits simple singularities. Assume that $K_{\cal F}+\Delta$ is pseudo-effective.

Then $(\cal F, \Delta)$ admits a minimal model $\psi\colon X\dashrightarrow Y$. Moreover, if $\cal G\coloneqq \psi_*\cal F$ and $\Gamma\coloneqq \psi_*\Delta$, then the following hold:
\begin{enumerate}
\item $\cal G$ admits simple singularities;
\item $(\cal G,\Gamma)$ is log canonical; 
\item if $\Theta \ge 0$ is a $\mathbb Q$-divisor on $X$ with $\cal F$-invariant support such that $(X,\Delta+\Theta)$ is log canonical, then 
$(Y,\psi_*(\Delta+\Theta))$ is log canonical. 
\end{enumerate}
\end{proposition}
\begin{proof}
Lemma \ref{l_simple_sing_are_quot} implies that $X$ has quotient singularities. 
In particular, $X$ is klt and $\mathbb Q$-factorial. 

If $K_{\cal F}+\Delta$ is nef then there is nothing to prove,  so we may assume that $K_{\cal F}+\Delta$ is not nef.
Let $R$ be a $(K_{\cal F}+\Delta)$-negative extremal ray.  By Theorem \ref{t_cone2} and Remark \ref{r_cone}, we may find
an $\cal F$-invariant curve $C$ spanning $R$. In particular, $C$ is a log canonical centre for $\cal F$. 
Since $(\cal F, \Delta)$ is log canonical, it follows
that no component of $\Delta$ is $\cal F$-invariant and  $\Delta\cdot C \geq 0$. Thus, 
$K_{\cal F}\cdot C<0$.

We may therefore apply Theorem \ref{t_div_contraction} and Theorem \ref{thm_flips_exist} to conclude that the contraction associated to $R$ exists and,
if the contraction is small, that the flip exists.  Call this step of the MMP $\phi\colon X\dashrightarrow X'$ and let $\cal F'$ be the induced foliation on $X'$.   Theorem \ref{t_div_contraction} and Theorem \ref{thm_flips_exist} (and their proofs) imply that 
 $X'$ is projective, 
 $\cal F'$ has simple singularities and that if $\Theta \ge 0$ is a $\mathbb Q$-divisor on $X$ with $\cal F$-invariant support such that $(X,\Delta+\Theta)$ is log canonical, then 
$(X',\phi_*(\Delta+\Theta))$ is log canonical. 
Moreover, Lemma \ref{l_negativity} implies  that $(\cal F', \Delta')$ is log canonical.  Thus, replacing
$X, \Delta$ and $\Theta$  by $X', \phi_*\Delta$ and $\phi_*\Theta$, we may continue this process.

Each divisorial contraction drops the Picard number by one, and so we can only contract a divisor
finitely many times. By Theorem \ref{thm_term} we can only have finitely many flips and so
this process must eventually terminate in our desired minimal model. 
\end{proof}

\begin{remark}
\label{rem_rel_mmp}
Let $p\colon X \to Z$ be a fibration between normal projective varieties. Let $(\cal F,\Delta)$ be a rank one foliated pair on $X$ with log canonical singularities and such that $\cal F$ admits simple singularities.

Suppose that $K_{\mathcal F}+\Delta$ is pseudo-effective over $Z$.  We can run a relative $(K_{\cal F}+\Delta)$-MMP over $Z$, call it $\psi\colon X \dashrightarrow Y/Z$ which terminates in a model where $K_{\psi_*\mathcal F}+\psi_*\Delta$ is nef over $Z$.   Indeed, the proof of Proposition \ref{prop_admissible_mmp} can be adapted to this setting by requring that at each step of the MMP we only contract/flip extremal rays which are $p^*H$-trivial, where $H$ is an ample divisor on $Z$.
\end{remark}

\subsection{Foliated plt blow ups}

In this section, we explain how to perform a foliated analogue of the classical plt blow up. We begin with the following:

\begin{lemma}\label{l_nd}
Let $X$ be a normal projective
threefold and let 
$(\cal F, \Delta)$ be a foliated pair on $X$ with log canonical singularities.
Let $E$ be a valuation which is exceptional over $X$ and such that  $a(E,\cal F,\Delta)<0$. 

Then $a(E,\cal F,\Delta)=-1$. 
In particular, if $a(E,\cal F,\Delta)>-1$   for any exceptional divisor $E$ over $X$ then
$(\cal F,\Delta)$ is canonical.
\end{lemma}

\begin{proof}
Let $p\colon Y \rightarrow X$ be the birational morphism whose existence is guaranteed by Theorem \ref{t_resolution}
 and such that $E$ is a divisor on $Y$. 
Let 
 $\cal F_Y\coloneqq p^{-1}\cal F$ and let $\Delta_Y \coloneqq p^{-1}_*\Delta$.  
  We may write 
\[
K_{{\cal F}_Y}+\Delta_Y +F' = p^*(K_{\cal F}+\Delta) +F''
\]
where $F', F''\ge 0$ are $p$-exceptional $\mathbb Q$-divisor with no common components. 
 After possibly passing to a higher resolution, 
we may assume that
$(\cal F_Y, \Delta_Y+F)$ is log canonical (cf. \cite[pag. 282, Corollary]{mp13}), where  $F\coloneqq \sum \epsilon(F_i)F_i$ and the sum runs over all the prime $p$-exceptional divisors.

Assume by contradiction that $a(E,\cal F,\Delta)\in (-1,0)$. 
In particular, $E$ is contained in the support of $F'$. 
Since $(\cal F,\Delta)$ is log canonical, 
it follows that $E$ is not $\cal F'$-invariant. Let $\epsilon>0$ be a positive rational number such that $(\cal F_Y,\Delta_Y+F'+\epsilon E)$ is log canonical. 
By Proposition \ref{prop_admissible_mmp}, $(\cal F_Y,\Delta_Y+F'+\epsilon E)$ admits a minimal model 
$\phi\colon Y\dashrightarrow X'$ over $X$, which, in particular, contracts $E$, 
contradicting Item (1) of Lemma \ref{lem_lc_contraction}. 
\end{proof}

\begin{theorem}
\label{thm_exist_plt_mod}
Let $X$ be a normal projective
 threefold  and let
$(\cal F, \Delta = \sum a_iD_i)$ be a foliated pair on $X$ where $a_i \in [0, \epsilon(D_i)]$.
 
Then there exists a birational morphism
$\pi\colon X' \rightarrow X$ such
that, if  $\cal F' \coloneqq\pi^{-1}\cal F$ and $\Delta'=\pi^{-1}_*\Delta$,
and $\{E_i\}$ is the set of all $\pi$-exceptional divisors then 
\begin{enumerate}
\item $\cal F'$ has simple singularities; 

\item $(X', \sum E_i)$ 
is log canonical, where the sum is over  all the $\pi$-exceptional divisors; and

\item there exists  a $\pi$-exceptional $\mathbb Q$-divisor $E'\ge 0$ on $X'$ such that 
\[
K_{\cal F'}+\Delta' + \sum \epsilon(E_i)E_i+ E' = \pi^*(K_{\cal F}+\Delta)
\]
and $(\cal F', \Delta'+\sum \epsilon(E_i)E_i)$ is log canonical.

\end{enumerate}
Moreover if $(\cal F,\Delta)$ is log canonical but not canonical at the generic point of a subvariety $P$ of  $X$ then

\begin{enumerate}
\item[(4)] there exists a unique prime $\pi$-exceptional divisor $E_0$ on $X'$ which is not $\cal F'$-invariant and which is centred on $P$;   and
\item[(5)] no other $\pi$-exceptional divisor has centre $=P$.
% there exists a $\pi$-exceptional $\mathbb Q$-divisor $E_1\ge 0$ on $X'$ 
 %whose centre on $X$ is not contained in  $P$ and such that $E'=E_0+E_1$. 
\end{enumerate}
We  call the morphism $\pi$ a  {\bf foliated plt blow up} of $(\cal F,\Delta)$.
\end{theorem}
\begin{proof}
Let $p\colon Y \rightarrow X$ be the birational morphism whose existence is guaranteed by Theorem \ref{t_resolution}.
Let  $\cal F_Y\coloneqq p^{-1}\cal F$ and let $\Delta_Y \coloneqq p^{-1}_*\Delta$. We may write 
\[
K_{{\cal F}_Y}+\Delta_Y +\sum \epsilon(E_i)E_i+F' = p^*(K_{\cal F}+\Delta) +F''
\]
where $F', F''\ge 0$ are $p$-exceptional $\mathbb Q$-divisor with no common components and $\{E_i\}$ is the set of all $p$-exceptional divisors. 
After possibly passing to a higher resolution, 
we may assume that
$(\cal F_Y, \Delta_Y+\sum \epsilon(E_i)E_i)$ is log canonical and that $(Y, \Delta_Y+\sum E_i)$ is log canonical  (cf. \cite[pag. 282, Corollary]{mp13}).

If   $(\cal F,\Delta)$ is log canonical but not canonical  at the generic point of a subvariety $P$ of $X$,
Lemma \ref{l_nd} implies that there exists an exceptional divisor $E_i$
centred over $P$ such that $\epsilon(E_i) = 1$ and $E_i$ is not contained 
in the support of $F'+F''$.

By Proposition \ref{prop_admissible_mmp} (see also Remark \ref{rem_rel_mmp}), we may run a $(K_{\cal F_Y}+\Delta_Y +  \sum \max\{\epsilon(E_i)-t, 0\}E_i)$-MMP over $X$
for any $t>0$ sufficiently small.  Let $\phi\colon Y\dashrightarrow X'$ be the output of this MMP. Let $\cal F'\coloneqq \phi_*\cal F_Y$ and let $E'\coloneqq \phi_* F'$. 
By Proposition \ref{prop_admissible_mmp}, we see that $(X',\sum \phi_*E_i)$ is log canonical.  
It is easy to verify that $X'$ and $\cal F'$ satisfy (1)-(5).
\end{proof}

This has the following useful consequence which allows us to reduce the MMP with log canonical singularities
to the MMP with canonical singularities.

\begin{corollary}
\label{cor_lc_fano}
Let $X$ be a projective threefold with log canonical singularities 
and let $\cal F$ be a foliation on $X$ with log canonical singularities.
Let $R$ be a $K_{\cal F}$-negative extremal ray and let $C$ be an $\cal F$-invariant curve such that $[C] \in R$.
Suppose that there exists a closed point $P \in C$ such that $\cal F$ is not canonical at $P$.

Then $\loc R = X$ and $R$ is $K_X$-negative.
\end{corollary}
\begin{proof}
Since $C$ is $K_{\cal F}$-negative,  it is not contained in $\sing {\mathcal F}$, see \cite[Fact II.d.3]{mcq04}. 
Proposition \ref{prop_adjunction} implies
that $\cal F$ is terminal at all points of $C \setminus P$.
Let $\pi\colon X' \rightarrow X$ be a foliated plt blow up of $\cal F$, whose existence is guaranteed by Theorem \ref{thm_exist_plt_mod},
and write $K_{\cal F'}+E = \pi^*K_{\cal F}$ where $E \geq 0$ and $\cal F' = \pi^{-1}\cal F$. In particular, $\mu_{E_0} E=1$ where $E_0$ is the unique $p$-exceptional divisor $E_0$ centred at $P$ and which is not $\cal F'$-invariant. 
By Lemma \ref{l_nondicritical} and since $\cal F$ is log canonical, it follows that no component of $E$ is centred on $C$. Since 
$\cal F$ is terminal at all points of $C \setminus P$, it follows that $E=E_0$. 

Then $K_{\cal F'}$ is not nef and there exists a curve $C'$ in $X'$ spanning a $K_{\cal F'}$-negative rational 
curve and such that $\pi(C')=C$. 

Notice that $K_{\cal F'}\cdot C' <0$.  Let $P' = E_0 \cap C'$.
Next, observe that $\cal F'$ has simple singularities and, therefore,  Lemma \ref{l_nondicritical} implies that for any exceptional divisor $E_1$ centred at a closed point of $E$, we have 
\[
a(E_1,\cal F')> a(E_1,\cal F)\ge \epsilon(E_1)=0.
\] 
Thus, 
$\cal F'$ is terminal at all closed points of $E_0$. In particular, $\cal F'$ is terminal
at $P'$, and so $\cal F'$ is terminal at all points of $C'$.

By Proposition \ref{prop_reeb_stability}, it follows that $C'$ moves in a family of pairwise disjoint curves covering $X'$.
Let $B$ be a general curve in such a family. Then 
\[
K_{X'}\cdot B = K_{\cal F'}\cdot B = -2.\]

We may write $K_{X'}+F = \pi^*K_X$ where $F$ is an exceptional $\pi$-divisor.  Since $X$ is log canonical, it follows that $\mu_{E_0}F\le 1$ and since $B\cdot G=0$ for every $\pi$-exceptional divisor $G$ which is $\cal F'$-invariant, it follows that 
\[
K_X\cdot \pi(B) = (K_{X'}+F)\cdot B \leq (K_{\cal F'}+E_0)\cdot B <0.\]
 Since $\pi(B)$ spans $R$, our result follows.
\end{proof}

We now show that, in the case of dimension three, Theorem \ref{t_cone2} holds without any $\bb Q$-factoriality hypothesis.

\begin{theorem}
\label{t_cone3}
Let $X$ be a normal projective threefold and let $(\mathcal F, \Delta)$ be a rank one 
foliated pair on $X$. 

Then there are $\mathcal F$-invariant rational curves $C_1,C_2,\dots$ not contained in $\sing {\cal F}$ 
such that 
\[
0<-(K_{\cal F}+\Delta)\cdot C_i\le 2\dim X\]
and
$$\overline{\rm NE}(X)=\overline{\rm NE}(X)_{K_{\cal F}+\Delta\ge 0}+Z_{-\infty}+
\sum_i \mathbb R_+[C_i]$$
where $Z_{-\infty}\subset \overline{\rm NE}(X)$ is a subset contained in the span of the images of
$\overline{\rm NE}(W) \rightarrow \overline{\rm NE}(X)$ where $W \subset X$
are the non-log canonical centres of $(\cal F, \Delta)$.
\end{theorem}

\begin{proof}
We use the notation of Theorem \ref{t_cone2} and its proof.
Let $p\colon X' \rightarrow X$ be a plt blow up of $(\cal F, \Delta)$, whose existence is guaranteed by Theorem \ref{thm_exist_plt_mod}, 
and write $K_{\cal F'}+\Delta' = p^*(K_{\cal F}+\Delta)$.
Notice that for any $(K_{\cal F}+\Delta)$-negative extremal ray $R$ there exists
a $(K_{\cal F'}+\Delta')$-negative extremal ray $R'$ with $p_*R' = R$.  Therefore,
we see that Theorem \ref{t_cone2} on $X'$ implies Theorem \ref{t_cone2} on $X$.
\end{proof}

\subsection{MMP with log canonical singularities}

We make note of an easy consequence of the negativity lemma which will 
nevertheless be crucial.

\begin{lemma}
\label{lem_neg_lem2}
Let $X$ be a projective variety and let $(\cal F, \Delta)$ be a rank one foliated 
pair with log canonical singularities.
Let $\phi\colon X \dashrightarrow X^+$ 
be a step of a $(K_{\cal F}+\Delta)$-MMP
and let $D \subset X$ be an $\cal F$-invariant divisor such that $\phi$
is an isomorphism at the generic point of $D$ and write $D^+ \coloneqq \phi_*D$. Let $\cal F^+$ be the foliation induced on $X^+$ and let $\Delta^+\coloneqq \phi_*\Delta$.  Write
\[(K_{\cal F}+\Delta)\vert_D = K_{\cal G}+\Theta\]
and
\[(K_{\cal F^+}+\Delta^+)\vert_{D^+} = K_{\cal G^+}+\Theta^+\]
where $(\cal G,\Theta)$ and $(\cal G^+,\Theta^+)$ are the induced foliated  pairs 
on $D$ and $D^+$, respectively. Let $W\xrightarrow{g}D$ and $W\xrightarrow{h}D^+$ be a resolution of $D \dashrightarrow D^+$.

Then $g^*(K_{\cal G}+\Theta)-h^*(K_{\cal G^+}+\Theta^+) \geq 0$ and is non-zero
if $\phi$ is not an isomorphism in a neighborhood of $D$.

In particular, the following hold:

\begin{enumerate}
\item If $K_{\cal G}+\Theta$ is not pseudo-effective then $K_{\cal G^+}+\Theta^+$
is not pseudoeffective.

\item If $K_{\cal G}+\Theta \equiv 0$ and $\phi$ is not an ismorphism
in a neighborhood of $D$ then $K_{\cal G^+}+\Theta^+$ is not pseudo-effective.
\end{enumerate}
\end{lemma}
\begin{proof}
The result follows immediately from the fact that $\phi$ is $(K_{\mathcal F}+\Delta)$-negative
and Proposition \ref{prop_adjunction}.
\end{proof}

\begin{theorem}\label{thm_flips_exist2}
Let $X$ be a $\mathbb Q$-factorial klt projective threefold
and let $\cal F$ be a rank one foliation on $X$ with  canonical  singularities. 
Let $R$ be a $K_{\cal F}$-negative extremal ray such that 
 $D\coloneqq \loc R $ has  dimension one.

Then the flipping contraction  $c_R\colon X\to Z$ associated to $R$  exists in the category of projective varieties.  Moreover, the flip $\phi\colon X \dashrightarrow X^+$ associated to $R$ exists and if $\cal F^+$ is the foliation induced on $X^+$ then 

\begin{enumerate}
\item $X^+$ is projective and has klt singularities;

\item $\cal F^+$ has canonical singularities and 
$\cal F^+$ is terminal at every point of $\exc \phi^{-1}$; and

\item $\rho(X/Z) = \rho(X^+/Z)= 1$.
\end{enumerate}
\end{theorem}
\begin{proof}
 Let
$C$ be a connected component of $\loc R$.
By Theorem \ref{t_cone2} and Remark \ref{r_cone}, we may assume that no component of $C$ is contained in $\sing \cal F$.
 By Lemma \ref{l_alg_flip_cont}, 
the contraction $f\colon X \rightarrow Z$ associated to $R$ exists in the category of algebraic spaces.

By Proposition \ref{prop_reeb_stability} and  Proposition \ref{prop_adjunction}, 
there exists a unique closed point $P\in C$ around which $\cal F$ is not terminal and  every irreducible component of $C$ passes through $P$. 
Let $p\colon Y \rightarrow X$ be a foliated plt blow up, whose existence is guaranteed by Theorem \ref{thm_exist_plt_mod},  let $\cal G \coloneqq p^{-1}\cal F$
and write 
\[
\exc p = \sum E_\ell+\sum F_j+\sum G_k
\]
where $p(E_\ell) = P$, $p(G_k)$ is an irreducible component of $C$ and $F_j$ are all the other exceptional divisors
which do not satisfy either of the previous conditions. Note that, by definition of a plt blow-up, every $p$-exceptional divisor maps to a canonical centre. Thus, since  $P$ is the only closed point in $C$ around which $\cal F$ is not terminal, it follows that the centre of $F_j$ is not contained in $C$. 

Since $\cal F$ admits canonical singularities, we have that $K_{\cal G}=p^*K_{\cal F}$ and  Lemma \ref{l_nondicritical} implies that $\exc p$ is $\cal G$-invariant. 
It follows that $K_{\cal G}\vert_{G_k}$ is not pseudoeffective for all $k$,
that $K_{\cal G}\vert_{E_\ell} \equiv 0$ for all $\ell$ and $K_{\cal G}\vert_{F_j}$
is numerically trivial over $X$ for all $j$.

By Proposition \ref{prop_admissible_mmp}, we may run a $K_{\cal G}$-MMP which only contracts/flips curves which are trivial with respect to $p^*H_R$.
This MMP will therefore be an MMP over $Z$, denote it by $\psi\colon Y \dashrightarrow Y^+$.
We observe the following facts:
\begin{itemize}
\item $\psi$ is  an isomorphism in a neighbourhood of a general fibre 
of the induced morphism $F_j \rightarrow p(F_j)$. 

\item $\psi$ contracts all the divisors $G_k$. 
Indeed, by Lemma \ref{lem_neg_lem2} if $Y_i\dashrightarrow Y_{i+1}$ is some intermediate 
step of the MMP, $\cal G_i$ is the induced foliation on $Y_i$ and $G^i_k \neq 0$ is the strict transform of $G_k$ on $Y_i$ then
$K_{\cal G_i}\vert_{G^i_k}$ is not pseudoeffective and so $\psi$ must
eventually contract $G_k$.

\item $\psi$ contracts all the $E_\ell$.  Indeed, again by Lemma \ref{lem_neg_lem2}, if $Y_i \dashrightarrow Y_{i+1}$
is some intermediate step of the MMP,   $\cal G_i$ is the induced foliation on $Y_i$ and $E^i_\ell \neq 0$ is the strict transform
of $E_\ell$ on $Y_i$ then either $Y \dashrightarrow Y_i$ is an isomorphism
in a neighbourhood of $E_\ell$, in which case $K_{\cal G_i}\vert_{E^i_\ell} \equiv 0$, 
or $Y \dashrightarrow Y_i$ is not an isomorphism near $E_\ell$.
 In the latter case,  if we choose $i$ to be the smallest positive integer such that $Y_0\coloneqq Y \dashrightarrow Y_{i}$ is not an isomorphism near $E_\ell$, 
then it follows that  $K_{\cal G_{i}}\vert_{E^{i}_\ell}$ is not pseudo-effective and arguing as in (2), we see
that $\psi$ contracts $E_\ell$.  Thus, our claim 
follows if we can show that for all $\ell$ there exists an $i_\ell$
such that $Y \dashrightarrow Y_{i_\ell}$ is not an isomorphism
near $E_\ell$.  This, however, follows from the fact that each connected component
of $\sum E_\ell$ has non-empty intersection  either with one of the divisor $G_k$ or with every irreducible
component in $p^{-1}(C)$ which is a curve dominating an irreducible component of $C$. 
Our claim then follows by proceeding by induction on the number of divisors $E_\ell$.  
\end{itemize}

Next, write $K_Y = \pi^*K_X +\sum a_jF_j +H$ where $H$ is supported on the $E_\ell$ and $G_k$.
Since $X$ is klt we may find an $\epsilon>0$ such that $a_j > -(1-\epsilon)$ for all $j$.
Let $F^+_j = \psi_*F_j$ and notice that $F^+_j \neq 0$ for all $j$.  Observe that
we still have morphisms $F^+_j \rightarrow p(F_j)$ and that $K_{\cal G^+}\vert_{F^+_j}$
is numerically trivial over the generic point of $p(F_j)$.

By the last property in Theorems \ref{t_div_contraction} and \ref{thm_flips_exist}
we know that $(Y^+, \sum F^+_j)$ is log canonical.
We may therefore run a $(K_{Y^+}+\sum F^+_j)$-MMP 
which only contracts/flips curves which are trivial with respect to
$K_{\cal G^+}$ and $\psi_*p^*H_R$, call this MMP
$\rho\colon Y^+ \dashrightarrow X^+$.  Observe that this will again be an MMP over $Z$
and that the following hold:

\begin{enumerate}
\item $\rho_*F^+_j = 0$ for all $j$, in particular, $f^+\colon X^+ \rightarrow Z$
is a small morphism.

\item Set $\cal F^+ = \rho_*\cal G^+$.  Then $K_{\cal F^+}$ is nef over $Z$.
\end{enumerate}

We claim that $f^+\colon X^+ \rightarrow Z$ is the desired flip.  
Let $\Sigma_1,\dots,\Sigma_\ell$ be the irreducible components of $\exc f^+$.

\begin{claim}
$[\Sigma_i]$ all span the same extremal ray $R^+ \subset \overline{\rm NE}(X^+)$.
\end{claim}
\begin{proof}[Proof of Claim]
Without loss of generality, we may assume that $K_X$ is ample over $Z$. Otherwise,
we would be able to realise the flipping contraction and flip as a consequence of the fact that 
$R$ is $(K_X+ D)$-negative 
for some suitable $\mathbb Q$-divisor $D$  such that $(X,  D)$ is klt.

Suppose for the sake of contradiction that the curves $\Sigma_1,\dots,\Sigma_\ell$ do not all span the same
extremal ray in $\overline{NE}(X^+)$.  Let $\rho\colon X^+ \dashrightarrow W$ be the birational contraction obtained by running 
 a $K_{X^+}$-MMP which only contracts/flips which are trivial with respect to  the strict transform of $H_R$. 
Observe that $X$ is  the log canonical model
of $W$ over $Z$, 
and so we have a morphism $W \rightarrow X$ which is small.  However,
$X$ is $\bb Q$-factorial and so $W \rightarrow X$ is necessarily an isomorphism.

We make the following general observation. 
Suppose that  $\phi\colon W_0 \dashrightarrow W_1$ is a $K_{W_0}$-flip which flips a curve $C_1$ and where $C_1^+$ is the flipped curve.
Suppose moreover there exists a curve $C_2 \subset W_0$ such that $C_2$ does not lie on $\bb R_+[C_1]$ and 
let $C^+_2 = \phi_*C_2$.  Then $C_2^+$ and $C_1^+$ do not lie on the same ray.  Indeed, let $M$ be a supporting
hyperplane to $\bb R_+[C_1]$ and let $M' = \phi_*M$.  Since $M$ is the pull back of a divisor
on the base of the flip we have that $0< M\cdot C_2=M'\cdot C_2^+$ and $0 = M\cdot C_1= M'\cdot C^+_1$, as required.

By inductively applying the above observation we see that
if $\Sigma^+_i$ denotes the strict transforms (resp. flipped curve)
of $\Sigma_i$, then 
not all the $\Sigma^+_i$ span the same  
ray in $\overline{NE}(X)$.
However, on the other hand, 
the $\Sigma^+_i$ are all $f$-exceptional and so all span $R$, a contradiction.
\end{proof}

Observe that the claim  implies that $K_{\cal F^+}$ is ample over $Z$.  Indeed, by construction
$K_{\cal F^+}$ is nef over $Z$ and it is necessarily not numerically trivial over $Z$
and so $K_{\cal F^+}\cdot \Sigma_i>0$ for all $i$ as required.

Next, observe that either $K_X$ is nef over $Z$
 or $-K_X$ is nef over $Z$.  If $-K_X$ is nef over $Z$ then, since $f$ is birational, it is also big over $Z$ and we may write
 $-K_X\sim_{\mathbb Q,f} A+E$ where $A$ is an ample $\mathbb Q$-divisor over $Z$ and $E\ge 0$. Thus, if $D:=\epsilon E$ for some sufficiently small rational number $\epsilon>0$, then $D \geq 0$, $-(K_X+D)$ is ample over $Z$ and $(X, D)$ is klt.  Thus, the contraction
of $R$ can be realised as a $(K_X+D)$-negative contraction, and so $Z$ is projective.
If $K_X$ is nef over $Z$ then $-K_{X^+}$ is nef over $Z$ and arguing as in the previous
case we may conclude that $Z$ is projective. In particular, $\rho(X/Z)=\rho(X^+/Z)=1$ and our claims follow. 
\end{proof}

\begin{theorem}\label{t_mmexist}
Let $X$ be a $\bb Q$-factorial 
projective threefold with klt singularities and let $(\cal F, \Delta)$ be a  log canonical foliated pair of rank one on $X$.
Assume that $K_{\cal F}+\Delta$ is pseudo-effective.

Then $(\cal F, \Delta)$ admits a minimal model.
\end{theorem}
\begin{proof}
If $K_{\cal F}+\Delta$ is nef there is nothing to show. So we may assume that $K_{\cal F}+\Delta$ is not nef.
Let $R$ be a $(K_{\cal F}+\Delta)$-negative extremal ray and let $H_R$ be a supporting hyperplane
to $R$. 
We want to show that the contraction, and possibly the flip, 
associated to $R$ exists.  Assuming this claim, we may argue as in Proposition \ref{prop_admissible_mmp} to
conclude that a minimal model exists.

Arguing as in Proposition \ref{prop_admissible_mmp},
we may again reduce to the case where we have a $\cal F$-invariant curve $C$ spanning $R$
which is $K_{\cal F}$-negative.
	By Theorem \ref{t_cone2}, we have that $C$ is not contained in $\sing \cal F$ and Proposition \ref{prop_adjunction} implies that there exists at most one closed point $P\in C$ at which $\cal F$ is singular.

Suppose that $\cal F$ has simple singularities in a neighbourhood of $C$. Then Theorem \ref{t_div_contraction} and Theorem \ref{thm_flips_exist} imply that the contraction, and possibly the flip,
of $R$ exists.

Now suppose that $\cal F$ is log canonical and not canonical at $P$.  In this case,  
Corollary \ref{cor_lc_fano} implies that $\loc R=X$, a contradiction.

Now suppose that $\cal F$ is canonical but not simple at $P$.
If $\loc R$ is a divisor, then Theorem \ref{t_div_contraction} implies the existence of a contraction.
Thus, we may assume that $\loc R$ is a curve and the claim follows from Theorem \ref{thm_flips_exist2}.
\end{proof}